# FAST RATES FOR SUPPORT VECTOR MACHINES USING GAUSSIAN KERNELS[1]

By Ingo Steinwart and Clint Scovel

*Los Alamos National Laboratory*

For binary classification we establish learning rates up to the order of $n^{-1}$ for support vector machines (SVMs) with hinge loss and Gaussian RBF kernels. These rates are in terms of two assumptions on the considered distributions: Tsybakov's noise assumption to establish a small estimation error, and a new geometric noise condition which is used to bound the approximation error. Unlike previously proposed concepts for bounding the approximation error, the geometric noise assumption does not employ any smoothness assumption.

**1. Introduction.** In recent years support vector machines (SVMs) have been the subject of many theoretical considerations. Despite this effort, their learning performance on restricted classes of distributions is still widely unknown. In particular, it is unknown under which nontrivial circumstances SVMs can guarantee *fast* learning rates. The aim of this work is to use concepts like Tsybakov's noise assumption and local Rademacher averages to establish learning rates up to the order of $n^{-1}$ for nontrivial distributions. In addition to these concepts that are used to deal with the stochastic part of the analysis we also introduce a geometric assumption for distributions that allows us to estimate the approximation properties of Gaussian RBF kernels. Unlike many other concepts introduced for bounding the approximation error, our geometric assumption is *not* in terms of smoothness but describes the concentration and the noisiness of the data-generating distribution near the decision boundary.

Let us formally introduce the statistical classification problem. To this end let us fix a subset $X \subset \mathbb{R}^d$. We write $Y := \{-1, 1\}$. Given a finite *training set*

Received December 2003; revised June 2006.

[1]Supported by the LDRD-ER program of the Los Alamos National Laboratory.

*AMS 2000 subject classifications.* Primary 68Q32; secondary 62G20, 62G99, 68T05, 68T10, 41A46, 41A99.

*Key words and phrases.* Support vector machines, classification, nonlinear discrimination, learning rates, noise assumption, Gaussian RBF kernels.

$T = ((x_1, y_1), \ldots, (x_n, y_n)) \in (X \times Y)^n$, the classification task is to predict the *label* $y$ of a new sample $(x, y)$. In the standard batch model it is assumed that the samples $(x_i, y_i)$ are i.i.d. according to an unknown (Borel) probability measure $P$ on $X \times Y$. Furthermore, the new sample $(x, y)$ is drawn from $P$ independently of $T$. Given a *classifier* $\mathcal{C}$ that assigns to every training set $T$ a measurable function $f_T : X \to \mathbb{R}$, the prediction of $\mathcal{C}$ for $y$ is $\operatorname{sign}_T f(x)$, where $\operatorname{sign}(0) := 1$. The quality of such a function $f$ is measured by the *classification risk*

$$\mathcal{R}_P(f) := P(\{(x, y) : \operatorname{sign} f(x) \neq y\}),$$

which should be as small as possible. The smallest achievable risk $\mathcal{R}_P := \inf\{\mathcal{R}_P(f) | f : X \to \mathbb{R} \text{ measurable}\}$ is called the *Bayes risk* of $P$ and a function attaining this risk is called a *Bayes decision function* and is denoted by $f_P$. Obviously, a good classifier should at least produce decision functions whose risks converge to the Bayes risk for *all* distributions $P$. This leads to the notion of *universally consistent* classifiers which is thoroughly treated in [14]. The next naturally arising question is whether there are classifiers which guarantee a specific convergence rate for *all* distributions. Unfortunately, this is impossible by a result of Devroye (see [14], Theorem 7.2). However, if one restricts consideration to certain *smaller* classes of distributions, such "learning rates," for example, in the form of

$$P^n(T \in (X \times Y)^n : \mathcal{R}_P(f_T) \leq \mathcal{R}_P + C(x) n^{-\beta}) \geq 1 - e^{-x}, \qquad n \geq 1, x \geq 1,$$

where $\beta > 0$ and $C(x) > 0$ are constants, exist for various classifiers. Typical assumptions for such classes of distributions are either in terms of the smoothness of the function $\eta(x) := P(y = 1 | x)$ (see, e.g., [19, 38]), or in terms of the smoothness of the "decision boundary" (see, e.g., [18, 35]). Moreover, the corresponding learning rates are slower than $n^{-1/2}$ if no additional assumptions on the *amount* of the noise in the labels, for example, on the distribution of the random variable

(1)             $$\min\{1 - \eta(x), \eta(x)\} = \tfrac{1}{2} - |\eta(x) - \tfrac{1}{2}|$$

around the critical level $1/2$, are imposed. On the other hand, [35] showed that ERM-type classifiers can learn faster than $n^{-1/2}$, if one quantifies how likely the noise in (1) is close to $1/2$ (see Definition 2.2 in the following section). Unfortunately, however the ERM classifier considered in [35] requires substantial knowledge on *how* to approximate the desired Bayes decision functions. Moreover, ERM classifiers are based on combinatorial optimization problems and hence they are usually hard to implement and in general there exist no efficient algorithms.

On the one hand SVMs do not share the implementation issues of ERM since they are based on a convex optimization (see, e.g., [12, 26] for algorithmic aspects). On the other hand, however, their *known* learning rates are rather unsatisfactory since either the assumptions on the distributions are too restrictive as in [28] or the established learning rates are too slow as in [37]. Our aim is to give SVMs a better theoretical foundation by establishing fast learning rates for a wide class of distributions. To this end we propose a geometric noise assumption (see Definition 2.3) which describes the concentration of the measure $|2\eta - 1| \, dP_X$—where $P_X$ is the marginal distribution of $P$ with respect to $X$—near the decision boundary. This assumption is then used to determine the approximation properties of Gaussian kernels which are used in the SVMs we consider. Provided that the tuning parameters are optimally chosen our main result then shows that the resulting learning rates for these classifiers can be as fast as $n^{-1}$.

The rest of this work is organized as follows: In Section 2 we introduce the main concepts of this work and then present our results. In Section 3 we recall some basic theory on reproducing kernel Hilbert spaces and prove a new covering number bound for Gaussian kernels that describes a trade-off between the kernel widths and the radii of the covering balls. In Section 4 we then show the approximation results that are related to our proposed geometric noise assumption. The last sections of the work contain the actual proof of our rates: In Section 5 we establish a general bound for ERM-type classifiers involving local Rademacher averages which is used to bound the estimation error in our analysis of SVMs. In order to apply this result we need "variance bounds" for SVMs which are established in Section 6. Interestingly, it turns out that sharp versions of these bounds depend on both Tsybakov's noise assumption and the approximation properties of the kernel used. Finally, we prove our learning rates in Section 7.

**2. Definitions and main results.** In this section we first recall some basic notions related to support vector machines which are needed throughout this text. In Section 2.2, we then present a covering number bound for Gaussian RBF kernels which will play an important role in our analysis of the estimation error of SVMs. In Section 2.3 we recall Tsybakov's noise assumption which will allow us to establish learning rates faster than $n^{-1/2}$. Then, in Section 2.4, we introduce the new geometric assumption that is used to estimate the approximation error for SVMs with Gaussian RBF kernels. Finally, we present and discuss our learning rates in Section 2.5.

2.1. *RKHSs, SVMs and basic definitions.* For two functions $f$ and $g$ we use the notation $f(\lambda) \preceq g(\lambda)$ to mean that there exists a constant $C > 0$ such that $f(\lambda) \leq Cg(\lambda)$ over some specified range of values of $\lambda$. We also use the notation $\succeq$ with similar meaning and the notation $\sim$ when both $\preceq$ and $\succeq$ hold. In particular, we use the same notation for sequences.

If not stated otherwise, $X$ always denotes a compact subset of $\mathbb{R}^d$ which is equipped with the Borel $\sigma$-algebra.

Recall (see, e.g., [1, 6]) that every positive definite kernel $k \colon X \times X \to \mathbb{R}$ has a unique reproducing kernel Hilbert space $H$ (RKHS) whose unit ball is denoted by $B_H$. Although we sometimes use generic kernels and RKHSs, we are mainly interested in Gaussian RBF kernels, which are the most widely used kernels in practice. Recall that these kernels are of the form

$$k_\sigma(x, x') = \exp(-\sigma^2 \|x - x'\|_2^2), \qquad x, x' \in X,$$

where $\sigma > 0$ is a free parameter whose *inverse* $1/\sigma$ is called the *width* of $k_\sigma$. We usually denote the corresponding RKHSs which are thoroughly described in [32] by $H_\sigma(X)$ or simply $H_\sigma$.

Let us now recall the definition of SVMs. To this end let $P$ be a distribution on $X \times Y$ and $l \colon Y \times \mathbb{R} \to [0, \infty)$ be the *hinge loss*, that is,

$$l(y, t) := \max\{0, 1 - yt\}, \qquad y \in Y, t \in \mathbb{R}.$$

Furthermore, we define the *l-risk* of a measurable function $f \colon X \to \mathbb{R}$ by

$$\mathcal{R}_{l,P}(f) := \mathbb{E}_{(x,y) \sim P} l(y, f(x)).$$

Now let $H$ be a RKHS over $X$ consisting of *measurable* functions. For $\lambda > 0$ we denote a solution of

$$(2) \qquad \underset{\substack{f \in H \\ b \in \mathbb{R}}}{\arg\min}(\lambda \|f\|_H^2 + \mathcal{R}_{l,P}(f + b))$$

by $(\tilde{f}_{P,\lambda}, \tilde{b}_{P,\lambda})$. Recall that $\tilde{f}_{P,\lambda}$ is uniquely determined (see, e.g., [30]), while in some situations this is not true for the *offset* $\tilde{b}_{P,\lambda}$. In general we thus assume that $\tilde{b}_{P,\lambda}$ is an arbitrary solution. However, for the (trivial) distributions that satisfy $P(\{y^*\}|x) = 1$ $P_X$-a.s. for some $y^* \in Y$ we explicitly set $\tilde{b}_{P,\lambda} := y^*$ in order to control the size of the offset. Furthermore, if $P$ is an empirical distribution with respect to a training set $T = ((x_1, y_1), \dots, (x_n, y_n))$ we write $\mathcal{R}_{l,T}(f)$ and $(\tilde{f}_{T,\lambda}, \tilde{b}_{T,\lambda})$. Note that in this case the above condition under which we set $\tilde{b}_{T,\lambda} := y^*$ means that all labels $y_i$ of $T$ are equal to $y^*$. An algorithm that constructs $(\tilde{f}_{T,\lambda}, \tilde{b}_{T,\lambda})$ for every training set $T$ is called an *SVM with offset*. Furthermore, for $\lambda > 0$ we denote the unique solution of

$$(3) \qquad \underset{f \in H}{\arg\min}(\lambda \|f\|_H^2 + \mathcal{R}_{l,P}(f))$$

by $f_{P,\lambda}$ and for empirical distributions based on a training set $T$ we again write $f_{T,\lambda}$. A corresponding algorithm is called an *SVM without offset*. Recall that under some assumptions on the RKHS used and the choice of the regularization parameter $\lambda$ it can be shown that both SVM variants are universally consistent (see [29, 31, 39]); however, no satisfying learning rates have been established yet.

We also emphasize that in many theoretical papers only SVMs without offset are considered since the offset often causes serious technical problems in the analysis. However, in practice usually SVMs with offset are used and therefore we feel that these algorithms should be considered in theory, too. As we will see, our techniques can be applied for both variants. The resulting rates coincide.

2.2. *Covering numbers for Gaussian RKHSs.* In order to bound the estimation error of SVMs we need a complexity measure for the RKHSs used, which is introduced in this section. To this end let $A \subset E$ be a subset of a Banach space $E$. The *covering numbers* of $A$ are defined by

$$\mathcal{N}(A, \varepsilon, E) := \min \left\{ n \geq 1 : \exists x_1, \ldots, x_n \in E \text{ with } A \subset \bigcup_{i=1}^n (x_i + \varepsilon B_E) \right\},$$

$\varepsilon > 0$, where $B_E$ denotes the closed unit ball of $E$. Moreover, for a bounded linear operator $S : E \to F$ between two Banach spaces $E$ and $F$, the covering numbers are $\mathcal{N}(S, \varepsilon) := \mathcal{N}(SB_E, \varepsilon, F)$.

Given a training set $T = ((x_1, y_1), \ldots, (x_n, y_n)) \in (X \times Y)^n$ we denote the space of all equivalence classes of functions $f : X \times Y \to \mathbb{R}$ with norm

$$(4) \qquad \|f\|_{L_2(T)} := \left( \frac{1}{n} \sum_{i=1}^n |f(x_i, y_i)|^2 \right)^{1/2}$$

by $L_2(T)$. In other words, $L_2(T)$ is an $L_2$-space with respect to the empirical measure of $T$. Note that for a function $f : X \times Y \to \mathbb{R}$ a canonical representative in $L_2(T)$ is its restriction $f_{|T}$. In addition, $L_2(T_X)$ denotes the space of all (equivalence classes of) square integrable functions with respect to the empirical measure of $x_1, \ldots, x_n$.

The proof of our learning rates uses the behavior of $\mathcal{N}(B_{H_\sigma(X)}, \varepsilon, L_2(T_X))$ in $\varepsilon$ and $\sigma$ in order to bound the estimation error. Unfortunately, all known results on covering numbers for Gaussian RBF kernels emphasize the role of $\varepsilon$ and hence we will establish in Section 3 the following result which describes a suitable trade-off between the influence of $\varepsilon$ and $\sigma$.

THEOREM 2.1. *Let $\sigma \geq 1$, $X \subset \mathbb{R}^d$ be a compact subset with nonempty interior, and $H_\sigma(X)$ be the RKHS of the Gaussian RBF kernel $k_\sigma$ on $X$. Then for all $0 < p \leq 2$ and all $\delta > 0$, there exists a constant $c_{p,\delta,d} > 0$ independent of $\sigma$ such that for all $\varepsilon > 0$ we have*

$$\sup_{T \in (X \times Y)^n} \log \mathcal{N}(B_{H_\sigma(X)}, \varepsilon, L_2(T_X)) \leq c_{p,\delta,d} \sigma^{(1-p/2)(1+\delta)d} \varepsilon^{-p}.$$

2.3. *Tsybakov's noise assumption.* Now we recall Tsybakov's noise condition, which describes the amount of noise in the labels. In order to motivate Tsybakov's assumption let us first observe that by equation (1) the function $|2\eta - 1|$ can be used to describe the noise in the labels of a distribution $P$. Indeed, in regions where this function is close to 1 there is only a small amount of noise, whereas function values close to 0 only occur in regions with a high level of noise. The following definition in which we use the convention $t^\infty := 0$ for $t \in (0, 1)$ describes the size of the latter regions:

DEFINITION 2.2.   Let $0 \leq q \leq \infty$ and $P$ be a probability measure on $X \times Y$. We say that $P$ has *Tsybakov noise exponent* $q$ if there exists a constant $C > 0$ such that for all sufficiently small $t > 0$ we have

$$(5) \qquad P_X(\{x \in X : |2\eta(x) - 1| \leq t\}) \leq C \cdot t^q.$$

Obviously, $P$ has Tsybakov noise exponent $q > 0$ if and only if $|2\eta - 1|^{-1} \in L_{q,\infty}(P_X)$, where $L_{q,\infty}$ denotes a Lorentz space (see [5]). It is also easy to see that $P$ has Tsybakov noise exponent $q'$ for all $q' < q$ if $P$ has Tsybakov noise exponent $q$. Furthermore, all distributions obviously have noise exponent 0. In the other extreme case $q = \infty$ the conditional probability $\eta$ is bounded away from $1/2$. In particular, noise-free distributions have exponent $q = \infty$. Furthermore, for $q < \infty$ it is easy to check that Definition 2.2 is satisfied if and only if (5) holds for *all* $t > 0$ and a possibly different constant $C$. Finally, note that (5) does not make any assumptions on the *location* of the noisy set, and hence we prefer the notion "noise condition" rather than the often used term "margin condition."

2.4. *A new geometric assumption for distributions.* In this section we introduce a condition for distributions that will allow us to estimate the approximation error for Gaussian RBF kernels. To this end let $l$ be the hinge loss function and $P$ be a distribution on $X$. Let

$$\mathcal{R}_{l,P} := \inf\{\mathcal{R}_{l,P}(f) | f : X \to \mathbb{R} \text{ measurable}\}$$

denote the smallest possible $l$-risk of $P$. Since functions achieving the minimal $l$-risk occur in many situations we indicate them by $f_{l,P}$ if no confusion regarding the nonuniqueness of this symbol can be expected. Furthermore, recall that $f_{l,P}$ has a shape similar to the Bayes decision function sign $f_P$ (see, e.g., [30]). Now, given a RKHS $H$ over $X$ we define the *approximation error function* with respect to $H$ and $P$ by

$$(6) \qquad a(\lambda) := \inf_{f \in H}(\lambda \|f\|_H^2 + \mathcal{R}_{l,P}(f) - \mathcal{R}_{l,P}), \qquad \lambda \geq 0.$$

Note that the obvious analogue of the approximation error function *with offset* is not greater than the above approximation error function *without offset* and hence we restrict our attention to the latter for simplicity.

For $\lambda > 0$, the approximation error function describes how well $\lambda \|f_{P,\lambda}\|_H^2 + \mathcal{R}_{l,P}(f_{P,\lambda})$ approximates $\mathcal{R}_{l,P}$. For example, it was shown in [31] that we have $\lim_{\lambda \to 0} a(\lambda) = 0$ for *all* $P$ if $X$ is a compact metric space and $H$ is dense in the space of continuous functions $C(X)$. However, in nontrivial situations there cannot exist a convergence rate which holds uniformly for all distributions $P$. Since $H_\sigma(X)$ is dense in $C(X)$ for compact $X \subset \mathbb{R}^d$ and all $\sigma > 0$ these statements are in particular true for the approximation error functions $a_\sigma(\cdot)$ of the Gaussian RBF kernels with *fixed* width $1/\sigma$. Moreover, we are not aware of any weak condition on $\eta$ or $P$ that ensures $a_\sigma(\lambda) \preceq \lambda^\beta$ for $\lambda \to 0$ and some $\beta > 0$, and the results of [27] indicate that such behavior of $a_\sigma(\cdot)$ may actually require very restrictive conditions. In the following we will therefore present a condition on $P$ that allows us to estimate $a_\sigma(\lambda)$ by $\lambda$ *and* $\sigma$. In particular it will turn out that $a_\sigma(\lambda) \to 0$ with a polynomial rate in $\lambda$ if we relate $\sigma$ to $\lambda$ in a certain manner. In order to introduce this assumption on $P$ we first define the classes of $P$ by $X_{-1} := \{x \in X : \eta(x) < \frac{1}{2}\}$, $X_1 := \{x \in X : \eta(x) > \frac{1}{2}\}$ and $X_0 := \{x \in X : \eta(x) = \frac{1}{2}\}$ for some choice of $\eta$. Now we define a distance function $x \mapsto \tau_x$ by

$$(7) \qquad \tau_x := \begin{cases} d(x, X_0 \cup X_1), & \text{if } x \in X_{-1}, \\ d(x, X_0 \cup X_{-1}), & \text{if } x \in X_1, \\ 0, & \text{otherwise,} \end{cases}$$

where $d(x, A)$ denotes the distance of $x$ to a set $A$ with respect to the Euclidean norm. Roughly speaking, $\tau_x$ measures the distance of $x$ to the "decision boundary." Now we can present the already announced geometric condition for distributions.

DEFINITION 2.3. Let $X \subset \mathbb{R}^d$ be compact and $P$ be a probability measure on $X \times Y$. We say that $P$ has *geometric noise exponent* $\alpha > 0$ if there exists a constant $C > 0$ such that

$$(8) \qquad \int_X |2\eta(x) - 1| \exp\left(-\frac{\tau_x^2}{t}\right) P_X(dx) \le C t^{\alpha d/2}, \qquad t > 0.$$

We say that $P$ has geometric noise exponent $\infty$ if it has geometric noise exponent $\alpha$ for all $\alpha > 0$.

Note that in the above definition we neither make any kind of smoothness assumption nor do we assume a condition on $P_X$ in terms of absolute continuity with respect to the Lebesgue measure. Instead, the integral condition (8) describes the concentration of the measure $|2\eta - 1| \, dP_X$ near the decision boundary in the sense that the less the measure is concentrated in this region the larger the geometric noise exponent can be chosen. The following example illustrates this.

EXAMPLE 2.4. Since $\exp(-t) \leq C_\alpha t^{-\alpha}$ holds for all $t > 0$ and a constant $C_\alpha > 0$ only depending on $\alpha > 0$, we easily see that (8) is satisfied whenever

$$(9) \qquad\qquad (x \mapsto \tau_x^{-1}) \in L_{\alpha d}(|2\eta - 1| \, dP_X),$$

where $L_{\alpha d}(|2\eta - 1| \, dP_X)$ denotes the usual Lebesgue space of functions that are $\alpha d$-integrable with respect to the measure $|2\eta - 1| \, dP_X$. Now, let us suppose $X_0 = \varnothing$ for a moment. In this case $\tau_x$ measures the distance to the class $x$ does not belong to. In particular, (9) holds for $\alpha = \infty$ if and only if the two classes $X_{-1}$ and $X_1$ have strictly positive distance. Moreover, if (9) holds for some $0 < \alpha < \infty$ the two classes may "touch," that is, the decision boundary $\partial X_{-1} \cap \partial X_1$ is nonempty. Consequently, we can easily construct distributions $P$ that have geometric noise exponent $\infty$ and touching classes, but also satisfy $f_P \notin H_\sigma(X)$ for all $\sigma > 0$. However, note that for such $P$ the measure $|2\eta - 1| \, dP_X$ must obviously have a very low concentration near the decision boundary.

We now describe a simple regularity condition on $\eta$ near the decision boundary that can be used to guarantee a geometric noise exponent.

DEFINITION 2.5. Let $X \subset \mathbb{R}^d$, $P$ be a distribution on $X \times Y$ and $\gamma > 0$. We say that $P$ has an *envelope of order* $\gamma$ if there is a constant $c_\gamma > 0$ such that for $P_X$-almost all $x \in X$ we have

$$(10) \qquad\qquad |2\eta(x) - 1| \leq c_\gamma \tau_x^\gamma.$$

Obviously, if $P$ has an envelope of order $\gamma$ then the graph of $x \mapsto 2\eta(x) - 1$ lies in a multiple of the envelope defined by $\tau_x^\gamma$ at the top and by $-\tau_x^\gamma$ at the bottom. Consequently, $\eta$ can be very irregular away from the decision boundary but cannot be discontinuous when crossing it. The rate of convergence of $\eta(x) \to 1/2$ for $\tau_x \to 0$ is described by $\gamma$.

Interestingly, for distributions having both an envelope of order $\gamma$ and a Tsybakov noise exponent $q$ we can bound the geometric noise exponent, as the following theorem, which is proved in Section 4, shows.

THEOREM 2.6. *Let $X \subset \mathbb{R}^d$ be compact and $P$ be a distribution on $X \times Y$ that has an envelope of order $\gamma > 0$ and a Tsybakov noise exponent $q \in [0, \infty)$. Then $P$ has geometric noise exponent $(q + 1)\gamma d^{-1}$ if $q \geq 1$, and geometric noise exponent $\alpha$ for all $\alpha < (q + 1)\gamma d^{-1}$ otherwise.*

Now the main result of this subsection which is proved in Section 4 shows that for distributions having a nontrivial geometric noise exponent we can bound the approximation error function for Gaussian RBF kernels.

THEOREM 2.7. *Let $\sigma > 0$, $X$ be the closed unit ball of the Euclidean space $\mathbb{R}^d$ and $a_\sigma(\cdot)$ be the approximation error function with respect to $H_\sigma(X)$. Furthermore, let $P$ be a distribution on $X \times Y$ that has geometric noise exponent $0 < \alpha < \infty$ with constant $C$ in* (8). *Then there is a constant $c_d > 0$ depending only on the dimension $d$ such that for all $\lambda > 0$ we have*

$$(11) \qquad a_\sigma(\lambda) \le c_d(\sigma^d \lambda + C(2d)^{\alpha d/2} \sigma^{-\alpha d}).$$

In order to let the right-hand side of (11) converge to zero it is necessary to assume both $\lambda \to 0$ and $\sigma \to \infty$. An easy consideration shows that the fastest convergence rate is achieved if $\sigma(\lambda) := \lambda^{-1/((\alpha+1)d)}$. In this case we have $a_{\sigma(\lambda)}(\lambda) \preceq \lambda^{\alpha/(\alpha+1)}$. In particular, we can obtain rates up to linear order in $\lambda$ for sufficiently benign distributions. The price for this good approximation property is, however, an increasing complexity of the hypothesis class $B_{H_{\sigma(\lambda)}}$, as we have seen in Theorem 2.1.

2.5. *Learning rates for SVMs using Gaussian RBF kernels.* With the help of the geometric noise assumption we can now present our learning rates for SVMs using Gaussian RBF kernels. Note again that these polynomial rates do *not* require a smoothness assumption on $P$. Furthermore note that we use the convention $\frac{a\infty+b}{c\infty+d} := \frac{a}{c}$ for $a, c \in (0, \infty)$, $b, d \in [0, \infty)$ in order to make the presentation compact.

THEOREM 2.8. *Let $X$ be the closed unit ball of $\mathbb{R}^d$, and $P$ be a distribution on $X \times Y$ with Tsybakov noise exponent $q \in [0, \infty]$ and geometric noise exponent $\alpha \in (0, \infty)$. We define*

$$\beta := \begin{cases} \dfrac{\alpha}{2\alpha+1}, & \text{if } \alpha \le \dfrac{q+2}{2q}, \\[2ex] \dfrac{2\alpha(q+1)}{2\alpha(q+2)+3q+4}, & \text{otherwise,} \end{cases}$$

*and $\lambda_n := n^{-(\alpha+1)/\alpha\beta}$ and $\sigma_n := n^{\beta/(\alpha d)}$ in both cases. Then for all $\varepsilon > 0$ there exists a $C > 0$ such that for all $x \ge 1$ and $n \ge 1$ the SVM without offset using the Gaussian RBF kernel $k_{\sigma_n}$ satisfies*

$$\Pr{}^*\big(T \in (X \times Y)^n : \mathcal{R}_P(f_{T,\lambda_n}) \le \mathcal{R}_P + Cx^2 n^{-\beta+\varepsilon}\big) \ge 1 - e^{-x},$$

*where $\Pr^*$ denotes the outer probability of $P^n$ in order to avoid measurability considerations. If $\alpha = \infty$ the latter inequality holds if $\sigma_n = \sigma$ is a constant with $\sigma > 2\sqrt{d}$. Finally, all results also hold for the SVM with offset.*

REMARK 2.9. *The above learning rates are faster than the "parametric" rate $n^{-1/2}$ if and only if $\alpha > (3q+4)/(2q)$. For $q = \infty$ the latter condition becomes $\alpha > 3/2$ and in an "intermediate" case $q = 1$ it becomes $\alpha > 7/2$.*

Remark 2.10. It is important to note that our techniques can also be used to establish rates for other definitions of the sequences $(\lambda_n)$ and $(\sigma_n)$. In fact, Theorem 2.7 guarantees $a_{\sigma_n}(\lambda_n) \to 0$ (which is necessary for our techniques to produce any rate) if $\sigma_n \to \infty$ and $\sigma_n^d \lambda_n \to 0$. In particular, if $\lambda_n := n^{-\iota}$ and $\sigma_n := n^\kappa$ for some $\iota, \kappa > 0$ with $\kappa d < \iota$, these conditions are satisfied and a conceptually easy but technically involved modification of our proof can produce rates for certain ranges of $\iota$ (and thus $\kappa$). In order to keep the presentation as short as possible we have omitted the details and focused on the best possible rates.

Remark 2.11. Unfortunately, the choice of $\lambda_n$ and $\sigma_n$ that yields the optimal rates within our techniques, requires to know the values of $\alpha$ and $q$, which are typically not available. Adaptive methods which do not require such knowledge are still unknown.

Remark 2.12. Theorem 2.7 and Theorem 2.8 establish results for *all* distributions having some geometric noise exponent. However, for certain distributions of this type the resulting rates are not satisfactory. For example consider the distribution $P$ on $X := [-1, 1]$ whose marginal distribution $P_X$ equals the uniform distribution and whose conditional distribution $\eta(x) := P(y = 1|x)$ satisfies $|2\eta(x) - 1| = |x|^\gamma$, $x \in X$, for some constant $\gamma \in (0, \infty)$. Then $P$ obviously has Tsybakov noise exponent $q := 1/\gamma$, and Theorem 2.6 or a simple modification of the proof of Theorem 2.7 shows that $P$ has geometric noise exponent $\alpha := 1 + \gamma$. Theorem 2.8 thus gives a rate of the form $n^{-\beta+\varepsilon}$ for $\beta = \frac{2q^2 + 4q + 2}{5q^2 + 10q + 4}$, which is never faster than $n^{-1/2}$. Though this is disappointing at first glance, it is not really surprising since the proof of Theorem 2.7 is not tailored to distributions having such simple decision functions. We believe that sharper bounds on the approximation error function (and thus faster learning rates) for this and other distributions are possible, but a detailed analysis is beyond the scope of this paper.

Remark 2.13. Another interesting but open question is whether the obtained rates are optimal for the class of considered distributions. In order to approach this question let us consider the case $\alpha = \infty$, which roughly speaking describes the case of almost no approximation error. In this case our rates are essentially of the form $n^{(q+1)/(q+2)}$, which coincides with the rates Tsybakov (see [35]) achieved for certain ERM classifiers based on hypothesis classes of small complexity. The latter rates in turn cannot be improved in a minimax sense for certain classes of distributions as was also shown in [35]. This discussion indicates that the techniques used for the stochastic part of our analysis may be strong enough to produce optimal results. However, if we consider the case $\alpha < \infty$ then the approximation error function described in Theorem 2.7 and its influence on the estimation error (see our proofs, in particular Section 5 and Section 7) have a significant impact on the obtained rates. Since the sharpness of Theorem 2.7 is unclear to us we make no conjecture regarding the optimality of our rates in the general case.

**3. Proof of Theorem 2.1.** The main goal of this section is to prove Theorem 2.1, which is done in Section 3.2. To this end we provide in Section 3.1 some RKHS theory which is used throughout this work.

3.1. *Some basic RKHS theory.* For the proofs of this section we have to recall some basic facts from the theory of RKHSs. To this end let $X \subset \mathbb{R}^d$ be a compact subset and $k : X \times X \to \mathbb{R}$ be a continuous and positive semi-definite kernel with RKHS $H$. Then $H$ consists of continuous functions on $X$ and for $f \in H$ we have $\|f\|_\infty \le K\|f\|_H$, where

$$(12) \qquad K := \sup_{x \in X} \sqrt{k(x, x)}.$$

Consequently, if the embedding of the RKHS $H$ into the space of continuous functions $C(X)$ is denoted by

$$(13) \qquad J_H : H \to C(X)$$

we have $\|J_H\| \le K$. Furthermore, let us recall the representation of $H$ based on Mercer's theorem (see [13]). To this end let $K_X : L_2(X) \to L_2(X)$ be the integral operator defined by

$$(14) \qquad K_X f(x) := \int_X k(x, x') f(x') \, dx', \qquad f \in L_2(X), x \in X,$$

where $L_2(X)$ denotes the $L_2$-space on $X$ with respect to the Lebesgue measure. Then it was shown in [13] that the unique square root $K_X^{1/2}$ of $K_X$ is an isometric isomorphism between $L_2(X)$ and $H$.

3.2. *Proof of Theorem 2.1.* In order to prove Theorem 2.1 we need the following result which bounds the covering numbers of $H_\sigma(X)$ with respect to $C(X)$.

THEOREM 3.1. *Let $\sigma \ge 1$, $0 < p < 2$ and $X \subset \mathbb{R}^d$ be a compact subset with nonempty interior. Then there is a constant $c_{p,d} > 0$ independent of $\sigma$ such that for all $\varepsilon > 0$ we have*

$$\log \mathcal{N}(B_{H_\sigma(X)}, \varepsilon, C(X)) \le c_{p,d}\sigma^{(1-p/4)d}\varepsilon^{-p}.$$

PROOF. Let $B_d$ be the closed unit ball of the Euclidean space $\mathbb{R}^d$ and $\mathring{B}_d$ be its interior. Then there exists an $r \ge 1$ such that $X \subset rB_d$. Now, it was recently shown in [32] that the restrictions $H_\sigma(rB_d) \to H_\sigma(X)$ and $H_\sigma(rB_d) \to H_\sigma(\mathring{B}_d)$ are both isometric isomorphisms. Consequently, in the following we assume without loss of generality that $X = B_d$ or $X = \mathring{B}_d$ and do not concern ourselves with the distinction of both cases.

Now let us write $H_\sigma := H_\sigma(X)$ and $J_\sigma := J_{H_\sigma} : H_\sigma \to C(X)$ in order to simplify notation. Furthermore, let $K_\sigma : L_2(X) \to L_2(X)$ be the integral operator of $k_\sigma$ defined as in (14), and $\|\cdot\|$ denote the norm in $L_2(X)$. According to [13], Theorem 3, page 27, for any $f \in H_\sigma$, we obtain

$$\inf_{\|K_\sigma^{-1}h\| \leq R} \|f - h\| \leq \frac{1}{R}\|K_\sigma^{-1/2}f\|^2 = \frac{1}{R}\|f\|_{H_\sigma}^2,$$

where we use the convention $\|K_\sigma^{-1}h\| = \infty$ if $h \notin K_\sigma L_2(X)$. Suppose now that $\mathcal{H} \subset L_2(X)$ is a dense Hilbert space with $\|h\| \leq \|h\|_\mathcal{H}$, and that we have $K_\sigma : L_2(X) \to \mathcal{H} \subset L_2(X)$ with $\|K_\sigma : L_2(X) \to \mathcal{H}\| \leq c_{\sigma,\mathcal{H}} < \infty$ for some constant $c_{\sigma,\mathcal{H}} > 0$. It follows that

$$\inf_{\|h\|_\mathcal{H} \leq c_{\sigma,\mathcal{H}}R} \|f - h\| \leq \inf_{\|K_\sigma^{-1}h\| \leq R} \|f - h\| \leq \frac{1}{R}\|f\|_{H_\sigma}$$

and hence

$$\inf_{\|h\|_\mathcal{H} \leq R} \|f - h\| \leq \frac{c_{\sigma,\mathcal{H}}}{R}\|f\|_{H_\sigma}^2.$$

By [27], Theorem 3.1 it follows that $f$ is contained in the real interpolation space $(L_2(X), \mathcal{H})_{1/2,\infty}$ (see [7] for the definition of an interpolation space) and its norm in this space satisfies $\|f\|_{1/2,\infty} \leq 2\sqrt{c_{\sigma,\mathcal{H}}}\|f\|_{H_\sigma}$. Therefore we obtain a continuous embedding

$$\Upsilon_1 : H_\sigma \to (L_2(X), \mathcal{H})_{1/2,\infty},$$

with $\|\Upsilon_1\| \leq 2\sqrt{c_{\sigma,\mathcal{H}}}$. If in addition a subset inclusion $(L_2(X), \mathcal{H})_{1/2,\infty} \subset C(X)$ exists which defines a continuous embedding

$$\Upsilon_2 : (L_2(X), \mathcal{H})_{1/2,\infty} \to C(X),$$

we have a factorization $J_\sigma = \Upsilon_2\Upsilon_1$ and can conclude

$$(15) \quad \log \mathcal{N}(B_{H_\sigma(X)}, \varepsilon, C(X)) = \log \mathcal{N}(J_\sigma, \varepsilon) \leq \log \mathcal{N}\left(\Upsilon_2, \frac{\varepsilon}{2\sqrt{c_{\sigma,\mathcal{H}}}}\right).$$

Consequently, to bound $\log \mathcal{N}(J_\sigma, \varepsilon)$ we need to select an $\mathcal{H}$, compute $c_{\sigma,\mathcal{H}}$ and bound $\log \mathcal{N}(\Upsilon_2, \varepsilon)$. To that end let $\mathcal{H} := W^m(\mathring{X})$ be the Sobolev space with norm

$$\|f\|_m^2 = \sum_{|\alpha| \leq m} \|D^\alpha f\|^2,$$

where $|\alpha| := \sum_{i=1}^{d} \alpha_i$, $D^\alpha := \prod_{i=1}^{d} \partial_i^{\alpha_i}$, and $\partial_i^{\alpha_i}$ denotes the $\alpha_i$th partial derivative in the $i$th coordinate of $\mathbb{R}^d$. By the Cauchy–Schwarz inequality we obtain

$$(16) \qquad \|D^\alpha K_\sigma f\|^2 \leq \|f\|^2 \int_X \int_X |D_x^\alpha k_\sigma(x, \acute{x})|^2 \, d\acute{x} \, dx,$$

where the notation $D_x^\alpha$ indicates that the differentiation takes place in the $x$ variable. To address the term $D_x^\alpha k_\sigma(x, \acute{x})$ we note that

$$D_x^\alpha(e^{-|x|^2}) = (-1)^{|\alpha|} e^{-|x|^2/2} h_\alpha(x),$$

where the multivariate Hermite functions $h_\alpha(x) = \prod_{i=1}^{d} h_{\alpha_i}(x_i)$ are products of the univariate functions. Since $\int_\mathbb{R} h_k^2(x) \, dx = 2^k k! \sqrt{\pi}$ (see, e.g., [11]) we obtain

$$(17) \qquad \begin{aligned} \int_{\mathbb{R}^d} |D_x^\alpha(e^{-|x|^2})|^2 \, dx &= \int_{\mathbb{R}^d} e^{-|x|^2} h_\alpha^2(x) \, dx \\ &\leq \int_{\mathbb{R}^d} h_\alpha^2(x) \, dx = 2^{|\alpha|} \alpha! \pi^{d/2}, \end{aligned}$$

where we have used the definition $\alpha! := \prod_{i=1}^{d} \alpha_i!$. Applying the translation invariance of $k_\sigma$, we obtain

$$\int_{\mathbb{R}^d} |D_x^\alpha k_\sigma(x, \acute{x})|^2 \, d\acute{x} = \int_{\mathbb{R}^d} |D_{\acute{x}}^\alpha k_\sigma(0, \acute{x})|^2 \, d\acute{x} = \int_{\mathbb{R}^d} |D_{\acute{x}}^\alpha(e^{-\sigma^2 |\acute{x}|^2})|^2 \, d\acute{x},$$

and by a change of variables we can apply inequality (17) to the integral on the right-hand side,

$$\int_{\mathbb{R}^d} |D_{\acute{x}}^\alpha(e^{-\sigma^2 |\acute{x}|^2})|^2 \, d\acute{x} = \sigma^{2|\alpha|-d} \int_{\mathbb{R}^d} |D_{\acute{x}}^\alpha(e^{-|\acute{x}|^2})|^2 \, d\acute{x} \leq \sigma^{2|\alpha|-d} 2^{|\alpha|} \alpha! \pi^{d/2}.$$

Hence we obtain

$$\int_X \int_X |D_x^\alpha k_\sigma(x, \acute{x})|^2 \, d\acute{x} \, dx \leq \theta(d) \sigma^{2|\alpha|-d} 2^{|\alpha|} \alpha! \pi^{d/2},$$

where $\theta(d)$ is the volume of $X$. Since $\sum_{|\alpha| \leq m} \alpha! \leq d^m m!^d$ and $\|K_\sigma f\|_m^2 = \sum_{|\alpha| \leq m} \|D^\alpha K_\sigma f\|^2$ we can therefore infer from (16) that for $\sigma \geq 1$ we have

$$(18) \qquad \|K_\sigma\| \leq \sqrt{\theta(d)} (2d)^{m/2} m!^{d/2} \sigma^{m-d/2} =: c_{\sigma, \mathcal{H}}.$$

Now let us consider $\Upsilon_2 : (L_2(X), W^m(\mathring{X}))_{1/2, \infty} \to C(X)$. According to Triebel [34], page 267, we have

$$(L_2(X), W^m(\mathring{X}))_{1/2, \infty} = (L_2(\mathring{X}), W^m(\mathring{X}))_{1/2, \infty} = B_{2, \infty}^{m/2}(\mathring{X})$$

isomorphically. Furthermore

$$(19) \qquad \log \mathcal{N}(B_{2, \infty}^{m/2}(\mathring{X}) \to C(X), \varepsilon) \leq c_{m, d} \varepsilon^{-2d/m}$$

for $m > d$ follows from a similar result of Birman and Solomjak ([8], cf. also [34]) for Slobodeckij (i.e., fractional Sobolev) spaces, where the constant $c_{m,d}$ depends only on $m$ and $d$. Consequently we obtain from (15), (18) and (19) that

$$\log \mathcal{N}(J_\sigma, \varepsilon) \le c_{m,d} \left( \frac{\varepsilon}{2\sqrt{c_{\sigma,\mathcal{H}}}} \right)^{-2d/m}$$

$$= c_{m,d} (4c_{\sigma,\mathcal{H}})^{d/m} \varepsilon^{-2d/m}$$

$$= \tilde{c}_{m,d} \sigma^{d - d^2/(2m)} \varepsilon^{-2d/m}$$

for all $m > d$ and new constants $\tilde{c}_{m,d}$ depending only on $m$ and $d$. Setting $m := 2d/p$ completes the proof of Theorem 3.1. $\square$

PROOF OF THEOREM 2.1. As before we write $H_\sigma := H_\sigma(X)$ and $J_\sigma := J_{H_\sigma} : H_\sigma \to C(X)$ in order to simplify notation. Furthermore recall for a training set $T \in (X \times Y)^n$ the space $L_2(T_X)$ introduced in Section 2.2. Now let $R_{T_X} : C(X) \to L_2(T_X)$ be the restriction map defined by $f \mapsto f_{|T_X}$. Obviously, we have $\|R_{T_X}\| \le 1$. Furthermore we define $I_\sigma := R_{T_X} \circ J_\sigma$ so that $I_\sigma : H_\sigma \to L_2(T_X)$ is the evaluation map. Then Theorem 3.1 and the product rule for covering numbers imply that

$$(20) \qquad \sup_{T \in Z^n} \log \mathcal{N}(I_\sigma, \varepsilon) \le c_{q,d} \sigma^{(1-q/4)d} \varepsilon^{-q}$$

for all $0 < q < 2$. To complete the proof of Theorem 2.1 we derive another bound on the covering numbers and interpolate the two. To that end observe that $I_\sigma : H_\sigma \to L_2(T_X)$ factors through $C(X)$ with both factors $J_s$ and $R_{T_X}$ having norm not greater than 1. Hence Proposition 17.3.7 in [23] implies that $I_\sigma$ is absolutely 2-summing with 2-summing norm not greater than 1. By König's theorem ([24], Lemma 2.7.2) we obtain for the approximation numbers $(a_k(I_\sigma))$ of $I_\sigma$ that $\sum_{k \ge 1} a_k^2(I_\sigma) \le 1$ for all $\sigma > 0$. Since the approximation numbers are decreasing it follows that $\sup_k \sqrt{k} a_k(I_\sigma) \le 1$. Using Carl's inequality between approximation and entropy numbers (see Theorem 3.1.1 in [10]) we thus find a constant $\tilde{c} > 0$ such that

$$(21) \qquad \sup_{T \in Z^n} \log \mathcal{N}(I_\sigma, \varepsilon) \le \tilde{c} \varepsilon^{-2}$$

for all $\varepsilon > 0$ and all $\sigma > 0$. Let us now interpolate the bound (21) with the bound (20). Since $\|I_\sigma : H_\sigma \to L_2(T_X)\| \le 1$ we only need to consider $0 < \varepsilon \le 1$. Let $0 < q < p < 2$ and $0 < a \le 1$. Then for $0 < \varepsilon < a$ we have

$$\log \mathcal{N}(I_\sigma, \varepsilon) \le c_{q,d} \sigma^{(1-q/4)d} \varepsilon^{-q} \le c_{q,d} \sigma^{(1-q/4)d} a^{p-q} \varepsilon^{-p},$$

and for $a \le \varepsilon \le 1$ we find

$$\log \mathcal{N}(I_\sigma, \varepsilon) \le \tilde{c} \varepsilon^{-2} \le \tilde{c} a^{p-2} \varepsilon^{-p}.$$

Since $\sigma \geq 1$ we can set $a := \sigma^{-((4-q)/(8-4q))d}$ and obtain

$$\log \mathcal{N}(I_\sigma, \varepsilon) \leq \tilde{c}_{q,d} \sigma^{(1-p/2)((8-2q)/(8-4q))d} \varepsilon^{-p},$$

where $\tilde{c}_{q,d}$ is a constant depending only on $q, d$. The proof is completed by choosing $q := \frac{4\delta}{1+2\delta}$ when $\delta < \frac{2p}{8-4p}$ and $q$ just smaller than $p$ otherwise. $\square$

**4. Proofs of Theorems 2.7 and 2.6.** In this section we prove Theorems 2.7 and 2.6, which both deal with the geometric noise exponent.

4.1. *Proof of Theorem 2.7.* Let us begin by recalling some facts about Gaussian RBF kernels. To this end let $H_\sigma(\mathbb{R}^d)$ be the RKHS of the Gaussian RBF kernel with parameter $\sigma$. Then it was shown in [32] that the linear operator $V_\sigma : L_2(\mathbb{R}^d) \to H_\sigma(\mathbb{R}^d)$ defined by

$$V_\sigma g(x) = \frac{(2\sigma)^{d/2}}{\pi^{d/4}} \int_{\mathbb{R}^d} e^{-2\sigma^2 \|x-y\|_2^2} g(y) \, dy, \qquad g \in L_2(\mathbb{R}^d), x \in \mathbb{R}^d,$$

is an isometric isomorphism. Consequently, we obtain

$$(22) \qquad a_\sigma(\lambda) = \inf_{g \in L_2(\mathbb{R}^d)} \lambda \|g\|_{L_2(\mathbb{R}^d)}^2 + \mathcal{R}_{l,P}(V_\sigma g) - \mathcal{R}_{l,P}, \qquad \lambda > 0.$$

In the following we will estimate the right-hand side of (22) by a judicious choice of $g$. To this end we need the following lemma, which in some sense enlarges the support of $P$ to ensure that all balls of the form $B(x, \tau_x)$ are contained in the (enlarged) support. This guarantee will then make it possible to control the behavior of $V_\sigma g$ by tails of spherical Gaussian distributions [see [28] for details].

LEMMA 4.1. *Let $X$ be a closed unit ball of $\mathbb{R}^d$ and $P$ be a probability measure on $X \times Y$ with regular conditional probability $\eta(x) = P(y = 1 | x)$, $x \in X$. On $\acute{X} := 3X$ we define*

$$(23) \qquad \acute{\eta}(x) = \begin{cases} \eta(x), & \text{if } |x| \leq 1, \\ \eta\left(\dfrac{x}{|x|}\right), & \text{otherwise.} \end{cases}$$

*We also write $\acute{X}_{-1} := \{x \in \acute{X} : \acute{\eta}(x) < \frac{1}{2}\}$ and $\acute{X}_1 := \{x \in \acute{X} : \acute{\eta}(x) > \frac{1}{2}\}$. Finally let $B(x, r)$ denote the open ball of radius $r$ about $x$ in $\mathbb{R}^d$. Then for $x \in X_1$ we have $B(x, \tau_x) \subset \acute{X}_1$ and for $x \in X_{-1}$ we have $B(x, \tau_x) \subset \acute{X}_{-1}$.*

PROOF. Let $x \in X_1$ and $x' \in B(x, \tau_x)$. If $x' \in X$ we have $|x - x'| < \tau_x$ which implies $\eta(x) > \frac{1}{2}$ by the definition of $\tau_x$. This shows $x' \in \acute{X}_1$. Now let us assume $|x'| > 1$. By $|\langle x, x' \rangle| \leq |x'|$ and Pythagoras' theorem we then obtain

$$\left| \frac{x'}{|x'|} - x \right|^2 \leq \left| x' - \frac{\langle x, x' \rangle x'}{|x'|^2} \right|^2 + \left| \frac{\langle x, x' \rangle x'}{|x'|^2} - x \right|^2 = |x' - x|^2.$$

Therefore, we have $|\frac{x'}{|x'|} - x| < \tau_x$, which implies $\acute{\eta}(x') = \eta(\frac{x'}{|x'|}) > \frac{1}{2}$.    □

Let us finally recall that Zhang showed in [39] that the hinge risk satisfies

$$(24) \qquad \mathcal{R}_{l,P}(f) - \mathcal{R}_{l,P} = \mathbb{E}_{P_X}(|2\eta - 1| \cdot |f - f_P|)$$

for all measurable $f : X \to [-1, 1]$. Now we are ready to prove Theorem 2.7.

PROOF OF THEOREM 2.7.  With the notation of Lemma 4.1 we fix a measurable $\acute{f}_P : \acute{X} \to [-1, 1]$ that satisfies $\acute{f}_P = 1$ on $\acute{X}_1$, $\acute{f}_P = -1$ on $\acute{X}_{-1}$ and $\acute{f}_P = 0$ otherwise. For $g := (\sigma^2/\pi)^{d/4} \acute{f}_P$ we then immediately obtain

$$(25) \qquad \|g\|_{L_2(\mathbb{R}^d)} \leq \left(\frac{81\sigma^2}{\pi}\right)^{d/4} \theta(d),$$

where $\theta(d)$ denotes the volume of $X$. Moreover, it is easy to see that $-1 \leq \acute{f}_P \leq 1$ implies $-1 \leq V_\sigma g \leq 1$. Since $P_X$ has support in $X$, (24) then yields

$$(26) \qquad \mathcal{R}_{l,P}(V_\sigma g) - \mathcal{R}_{l,P} = \mathbb{E}_{P_X}(|2\eta - 1| \cdot |V_\sigma g - f_P|).$$

In order to bound $|V_\sigma g(x) - f_P(x)|$ for $x \in X_1$ we observe

$$
\begin{aligned}
(27) \qquad V_\sigma g(x) &= \left(\frac{2\sigma^2}{\pi}\right)^{d/2} \int_{\mathbb{R}^d} e^{-2\sigma^2 \|x-y\|_2^2} \acute{f}_P(y) \, dy \\
&= \left(\frac{2\sigma^2}{\pi}\right)^{d/2} \int_{\mathbb{R}^d} e^{-2\sigma^2 \|x-y\|_2^2} (\acute{f}_P(y) + 1) \, dy - 1 \\
&\geq \left(\frac{2\sigma^2}{\pi}\right)^{d/2} \int_{B(x,\tau_x)} e^{-2\sigma^2 \|x-y\|_2^2} (\acute{f}_P(y) + 1) \, dy - 1.
\end{aligned}
$$

Now remember that Lemma 4.1 showed $B(x, \tau_x) \subset \acute{X}_1$ for all $x \in X_1$, so that (27) implies

$$
\begin{aligned}
(28) \qquad V_\sigma g(x) &\geq 2 \left(\frac{2\sigma^2}{\pi}\right)^{d/2} \int_{B(x,\tau_x)} e^{-2\sigma^2 \|x-y\|_2^2} \, dy - 1 \\
&= 1 - 2P_{\gamma_\sigma}(|u| \geq \tau_x),
\end{aligned}
$$

where $\gamma_\sigma = (2\sigma^2/\pi)^{d/2} e^{-2\sigma^2 |u|^2} \, du$ is a spherical Gaussian in $\mathbb{R}^d$. According to the tail bound [17], inequality (3.5) on page 59, we have $P_{\gamma_\sigma}(|u| \geq r) \leq 4e^{-\sigma^2 r^2/2d}$ and consequently we obtain

$$1 \geq V_\sigma g(x) \geq 1 - 8e^{-\sigma^2 \tau_x^2/2d}, \qquad x \in X_1.$$

Since for $x \in X_{-1}$ we can obtain an analogous estimate, we conclude

$$|V_\sigma g(x) - f_P(x)| \leq 8e^{-\sigma^2 \tau_x^2/2d}$$

for all $x \in X_1 \cup X_{-1}$. Consequently (26) and the geometric noise assumption for $t := \frac{2d}{\sigma^2}$ yield

$$(29) \qquad \begin{aligned} \mathcal{R}_{l,P}(V_\sigma g) - \mathcal{R}_{l,P} &\leq 8\mathbb{E}_{x \sim P_X}(|2\eta(x) - 1|e^{-\sigma^2 \tau_x^2 / 2d}) \\ &\leq 8C(2d)^{\alpha d/2}\sigma^{-\alpha d}, \end{aligned}$$

where $C$ is the constant in (8). Combining (29), (25) and (22) now yields the assertion. $\square$

4.2. *Proof of Theorem* 2.6. In this subsection, all Lebesgue and Lorentz spaces (see, e.g., [5]) and their norms are with respect to the measure $P_X$.

PROOF OF THEOREM 2.6. Let us first consider the case $q \geq 1$ where we can apply the Hölder inequality for Lorentz spaces [22], which states

$$\|fg\|_1 \leq \|f\|_{q,\infty}\|g\|_{q',1}$$

for all $f \in L_{q,\infty}$, $g \in L_{q',1}$ and $q'$ defined by $\frac{1}{q} + \frac{1}{q'} = 1$. Applying this inequality gives

$$(30) \qquad \begin{aligned} \mathbb{E}_{x \sim P_X}&(|2\eta(x) - 1|e^{-\tau_x^2/t}) \\ &\leq \|(2\eta - 1)^{-1}\|_{q,\infty}\|x \mapsto (2\eta(x) - 1)^2 e^{-\tau_x^2/t}\|_{q',1} \\ &\leq C\|(2\eta - 1)^2 e^{-(|2\eta - 1|/c_\gamma)^{2/\gamma} t^{-1}}\|_{q',1}, \end{aligned}$$

where in the last estimate we used the Tsybakov assumption (5) and the fact that $P$ has an envelope of order $\gamma$. Let us write $h(x) := |2\eta(x) - 1|^{-1}$, $x \in X$, and $b := t(c_\gamma)^{2/\gamma}$ so that

$$|2\eta(x) - 1|^2 e^{-(|2\eta - 1|/c_\gamma)^{2/\gamma} t^{-1}} = g(h(x)),$$

where $g(s) := s^{-2}e^{-(s^{-2/\gamma})/b}$ for all $s \geq 1$. Now it is easy to see that $g : [1, \infty) \to [0, \infty)$ is strictly increasing if $0 < b \leq \frac{2}{3\gamma}$, and hence we can extend $g$ to a strictly increasing, continuous and invertible function on $[0, \infty)$ in this case. Let such an extension also be denoted by $g$. Then for this extension we have

$$(31) \qquad P_X(g \circ h > \tau) = P_X(h > g^{-1}(\tau)).$$

Now for a function $f : X \to [0, \infty)$ recall the nonincreasing rearrangement

$$f^*(u) := \inf\{\sigma \geq 0 : P_X(f > \sigma) \leq u\}, \qquad u > 0,$$

of $f$ which can be used to define Lorentz norms (see, e.g., [5]). For $u > 0$ equation (31) then yields

$$(g \circ h)^*(u) = g(\inf\{g^{-1}(\sigma) : P_X(h > g^{-1}(\sigma)) \leq u\}) = g \circ h^*(u).$$

Now, inequality (5) implies $P_X(h \geq (\frac{C}{u})^{1/q}) \leq u$ for all $u > 0$. Therefore, we find

$$h^*(u) \leq \inf\{\sigma \geq 0 : P_X(h \geq \sigma) \leq u\} \leq \left(\frac{C}{u}\right)^{1/q}$$

for all $0 < u < 1$. Since $(g \circ h)^* = g \circ h^*$ and $g$ is increasing we hence have

$$(g \circ h)^*(u) \leq g\left(\left(\frac{C}{u}\right)^{1/q}\right)$$

for all $0 < u < 1$. Now, for fixed $\hat{\alpha} > 0$ the bound $e^{-x} \preceq \frac{x^{-\hat{\alpha}}}{\ln^2(x)+1}$ on $(0, \infty)$ implies

$$g(s) \preceq b^{\hat{\alpha}} \frac{s^{2(\hat{\alpha}/\gamma-1)}}{\ln^2\left(s^{-2/\gamma}b^{-1}\right)+1}$$

for $s \in [1, \infty)$. Using the fact that $(g \circ h)^*(u) = 0$ holds for all $u \geq 1$, we hence obtain

$$(g \circ h)^*(u) \preceq b^{\hat{\alpha}} \frac{u^{2/q(1-\hat{\alpha}/\gamma)}}{\ln^2\left((u/C)^{2/(q\gamma)}b^{-1}\right)+1}$$

for $u > 0$ if we assume without loss of generality that $C \geq 1$. Let us define $\hat{\alpha} := \gamma\frac{q+1}{2}$. Then we find $\frac{1}{q'} + \frac{2}{q}(1 - \frac{\hat{\alpha}}{\gamma}) = 0$ and consequently for $b \leq \frac{2}{3\gamma}$, that is, $t \leq \frac{2}{3\gamma(c_\gamma)^{2/\gamma}}$, we obtain

$$
\begin{aligned}
(32) \quad \|g \circ h\|_{q',1} &= \int_0^\infty u^{1/q'-1}(g \circ h)^*(u)\, du \\
&\preceq b^{\hat{\alpha}} \int_0^\infty \frac{u^{-1}}{\ln^2\left((u/C)^{2/(q\gamma)}b^{-1}\right)+1}\, du \preceq t^{\gamma(q+1)/2}
\end{aligned}
$$

by the definition of $b$. Since we also have $\mathbb{E}_{P_X}(|2\eta(x)-1|e^{-\tau_x^2/t}) \leq 1$ for all $t > 0$, estimate (30) together the definition of $g$ and (32) yields the assertion in the case $q \geq 1$.

Let us now consider the case $0 \leq q < 1$ where the Hölder inequality in Lorentz space cannot be used. Then for all $t, \tau \geq 0$ we have

$$
\begin{aligned}
(33) \quad \mathbb{E}_{x \sim P_X}&(|2\eta(x)-1|e^{-\tau_x^2/t}) \\
&= \int_{|2\eta-1| \leq \tau} |2\eta(x)-1|e^{-\tau_x^2/t}P_X(dx) \\
&\quad + \int_{|2\eta-1| > \tau} |2\eta(x)-1|e^{-\tau_x^2/t}P_X(dx) \\
&\leq C\tau^{q+1} + \exp\left(-\left(\frac{\tau}{c_\gamma}\right)^{2/\gamma} t^{-1}\right),
\end{aligned}
$$

where we have used the Tsybakov assumption (5) and the fact that $P$ has an envelope of order $\gamma$. Let us define $\tau$ by $\tau^{q+1} := \exp(-(\frac{\tau}{c_\gamma})^{2/\gamma} t^{-1})$. For $\hat{a} := (c_\gamma)^{2/\gamma}(q+1)$ and small $t$ this definition implies

$$\tau \leq \left(\frac{\hat{a}\gamma}{2}\right)^{\gamma/2} \left(t \ln \frac{1}{\hat{a}t}\right)^{\gamma/2},$$

and hence the assertion follows from (33) for the case $0 < q < 1$. $\square$

## 5. The estimation error of ERM-type classifiers.

To bound the estimation error in the proof of Theorem 2.8 we now establish a concentration inequality for ERM-type algorithms using a variant of Talagrand's concentration inequality together with local Rademacher averages (see, e.g., [2, 4, 21]). Our approach is inspired by [3]. However, due to the regularization term $\lambda\|f\|_H^2$ in the definition of SVMs we need a more general result than that of [3].

This section is organized as follows: In Section 5.1 we present the required modification of the result of [3]. Then in Section 5.2 we bound the resulting local Rademacher averages.

### 5.1. *Bounding the estimation error for ERM-type algorithms.*

We first have to introduce some notation. To this end let $\mathcal{F}$ be a class of bounded measurable functions from $Z$ to $\mathbb{R}$ such that $\mathcal{F}$ is *separable* with respect to $\|\cdot\|_\infty$. Given a probability measure $P$ on $Z$ we define the modulus of continuity of $\mathcal{F}$ by

$$\omega_n(\mathcal{F}, \varepsilon) := \omega_{P,n}(\mathcal{F}, \varepsilon) := \mathbb{E}_{T \sim P^n} \left( \sup_{\substack{f \in \mathcal{F}, \\ \mathbb{E}_P f^2 \leq \varepsilon}} |\mathbb{E}_P f - \mathbb{E}_T f| \right), \qquad \varepsilon > 0,$$

where we note that the supremum is, as a function from $Z$ to $\mathbb{R}$, measurable by the separability assumption on $\mathcal{F}$. Now, a function $L : \mathcal{F} \times Z \to [0, \infty)$ is called a *loss function* if $L \circ f := L(f, \cdot)$ is measurable for all $f \in \mathcal{F}$. Given a probability measure $P$ on $Z$ we indicate by $f_{P,\mathcal{F}} \in \mathcal{F}$ a minimizer of

$$f \mapsto \mathcal{R}_{L,P}(f) := \mathbb{E}_{z \sim P} L(f, z).$$

Throughout this paper $\mathcal{R}_{L,P}(f)$ is called the $L$-risk of $f$. If $P$ is an empirical measure with respect to $T \in Z^n$ we write $f_{T,\mathcal{F}}$ and $\mathcal{R}_{L,T}(\cdot)$ as usual. For simplicity, we assume throughout this section that $f_{P,\mathcal{F}}$ and $f_{T,\mathcal{F}}$ do exist. Furthermore, although there may be multiple solutions we use a single symbol for them whenever no confusion regarding the nonuniqueness of this symbol can be expected. An algorithm that produces solutions $f_{T,\mathcal{F}}$ is called an *empirical $L$-risk minimizer*. Moreover, if $\mathcal{F}$ is convex, we say that $L$ is convex if $L(\cdot, z)$ is convex for all $z \in Z$. Finally, $L$ is called *line-continuous* if for all $z \in Z$ and all $f, \hat{f} \in \mathcal{F}$ the function $t \mapsto L(tf + (1-t)\hat{f}, z)$ is continuous on $[0,1]$. If $\mathcal{F}$ is a vector space then every convex $L$ is line-continuous. Now the main result of this section reads as follows:

THEOREM 5.1. *Let $\mathcal{F}$ be a convex set of bounded measurable functions from $Z$ to $\mathbb{R}$, and let $L : \mathcal{F} \times Z \to [0, \infty)$ be a convex and line-continuous loss function. For a probability measure $P$ on $Z$ we define*

$$\mathcal{G} := \{L \circ f - L \circ f_{P, \mathcal{F}} : f \in \mathcal{F}\}.$$

*Suppose that there are constants $c \geq 0$, $0 < \alpha \leq 1$, $\delta \geq 0$ and $B > 0$ with $\mathbb{E}_P g^2 \leq c(\mathbb{E}_P g)^\alpha + \delta$ and $\|g\|_\infty \leq B$ for all $g \in \mathcal{G}$. Furthermore, assume that $\mathcal{G}$ is separable with respect to $\|\cdot\|_\infty$. Let $n \geq 1$, $x \geq 1$ and $\varepsilon > 0$ with*

$$(34) \qquad \varepsilon \geq 10 \max \left\{ \omega_n(\mathcal{G}, c\varepsilon^\alpha + \delta), \sqrt{\frac{\delta x}{n}}, \left(\frac{4cx}{n}\right)^{1/(2-\alpha)}, \frac{Bx}{n} \right\}.$$

*Then we have*

$$\Pr^*(T \in Z^n : \mathcal{R}_{L,P}(f_{T, \mathcal{F}}) < \mathcal{R}_{L,P}(f_{P, \mathcal{F}}) + \varepsilon) \geq 1 - e^{-x}.$$

REMARK 5.2. Theorem 5.1 has been proved in [3] for $\delta = 0$, where it was used to find learning rates faster than $n^{-1/2}$ for certain ERM-type algorithms. At first glance such fast rates are impossible if $\delta > 0$. However, we will see later that for SVMs we have $\delta = a_\sigma^\kappa(\lambda)$ for a suitable $\kappa > 0$ depending on both Tsybakov's and the geometric noise exponent, and hence we have $\delta \to 0$ for $n \to \infty$.

As already mentioned, the proof of Theorem 5.1 is based on Talagrand's concentration inequality in [33] and its refinements in [16, 20, 25]. The version below of this inequality is derived from Bousquet's result in [9] using a little trick presented in [2], Lemma 2.5.

THEOREM 5.3. *Let $P$ be a probability measure on $Z$ and $\mathcal{H}$ be a set of bounded measurable functions from $Z$ to $\mathbb{R}$ which is separable with respect to $\|\cdot\|_\infty$ and satisfies $\mathbb{E}_P h = 0$ for all $h \in \mathcal{H}$. Furthermore, let $b > 0$ and $\tau \geq 0$ be constants with $\|h\|_\infty \leq b$ and $\mathbb{E}_P h^2 \leq \tau$ for all $h \in \mathcal{H}$. Then for all $x \geq 1$ and all $n \geq 1$ we have*

$$P^n\left(T \in Z^n : \sup_{h \in \mathcal{H}} \mathbb{E}_T h > 3\mathbb{E}_{T' \sim P^n} \sup_{h \in \mathcal{H}} \mathbb{E}_{T'} h + \sqrt{\frac{2x\tau}{n}} + \frac{bx}{n}\right) \leq e^{-x}.$$

This concentration inequality is used to prove the following lemma which is a generalized version of Lemma 13 in [3].

LEMMA 5.4. *Let $P$ be a probability measure on $Z$ and $\mathcal{G}$ be a set of bounded measurable functions from $Z$ to $\mathbb{R}$ which is separable with respect to $\|\cdot\|_\infty$. Let $c \geq 0$, $0 < \alpha \leq 1$, $\delta \geq 0$ and $B > 0$ be constants with $\mathbb{E}_P g^2 \leq c(\mathbb{E}_P g)^\alpha + \delta$ and $\|g\|_\infty \leq B$ for all $g \in \mathcal{G}$. Furthermore, assume that for all $T \in Z^n$ and all $\varepsilon > 0$ for which for some $g \in \mathcal{G}$ we have*

$$\mathbb{E}_T g \leq \varepsilon/20 \quad and \quad \mathbb{E}_P g \geq \varepsilon$$

*there exists a $g^* \in \mathcal{G}$ which satisfies*

$$\mathbb{E}_T g^* \leq \varepsilon/20 \quad and \quad \mathbb{E}_P g^* = \varepsilon.$$

*Then for all $n \geq 1$, $x \geq 1$, and all $\varepsilon > 0$ satisfying (34), we have*

$$\Pr{}^*(T \in Z^n : \text{ for all } g \in \mathcal{G} \text{ with } \mathbb{E}_T g \leq \varepsilon/20 \text{ we have } \mathbb{E}_P g < \varepsilon) \geq 1 - e^{-x}.$$

PROOF. We define $\mathcal{H} := \{\mathbb{E}_P g - g : g \in \mathcal{G}, \mathbb{E}_P g = \varepsilon\}$. Obviously, we have $\mathbb{E}_P h = 0$, $\|h\|_\infty \leq 2B$, and $\mathbb{E}_P h^2 = \mathbb{E}_P g^2 - (\mathbb{E}_P g)^2 \leq c\varepsilon^\alpha + \delta$ for all $h \in \mathcal{H}$. Moreover, since it is also easy to verify that $\mathcal{H}$ is separable with respect to $\|\cdot\|_\infty$, our assumption on $\mathcal{G}$ yields

$$\Pr{}^*(T \in Z^n : \exists g \in \mathcal{G} \text{ with } \mathbb{E}_T g \leq \varepsilon/20 \text{ and } \mathbb{E}_P g \geq \varepsilon)$$

$$\leq \Pr{}^*(T \in Z^n : \exists g \in \mathcal{G} \text{ with } \mathbb{E}_P g - \mathbb{E}_T g \geq 19\varepsilon/20 \text{ and } \mathbb{E}_P g = \varepsilon)$$

$$\leq P^n\left(T \in Z^n : \sup_{h \in \mathcal{H}} \mathbb{E}_T h \geq 19\varepsilon/20\right).$$

Note that since $\mathcal{H}$ is separable with respect to $\|\cdot\|_\infty$, the set on the last line is actually measurable. In order to bound the last probability we will apply Theorem 5.3. To this end we have to show

$$\frac{19\varepsilon}{20} > 3\mathbb{E}_{T' \sim P^n} \sup_{h \in \mathcal{H}} \mathbb{E}_{T'} h + \sqrt{\frac{2x\tau}{n}} + \frac{bx}{n}.$$

Our assumptions on $\varepsilon$ imply

$$(35) \quad \varepsilon \geq 10\mathbb{E}_{T' \sim P^n}\left(\sup_{\substack{g \in \mathcal{G}, \\ \mathbb{E}_P g^2 \leq c\varepsilon^\alpha + \delta}} |\mathbb{E}_P g - \mathbb{E}_{T'} g|\right) \geq 10\mathbb{E}_{T' \sim P^n} \sup_{h \in \mathcal{H}} \mathbb{E}_{T'} h.$$

Furthermore, since $10 \geq (\frac{60}{19})^2$ and $0 < \alpha \leq 1$ we have

$$(36) \quad \varepsilon \geq 10\left(\frac{4cx}{n}\right)^{1/(2-\alpha)} \geq \left(\frac{60}{19}\right)^{2/(2-\alpha)}\left(\frac{4cx}{n}\right)^{1/(2-\alpha)}.$$

If $\delta \leq c\varepsilon^\alpha$ a simple calculation hence shows $\frac{19}{60}\varepsilon \geq \sqrt{\frac{2(c\varepsilon^\alpha + \delta)x}{n}}$. Furthermore, if $\delta > c\varepsilon^\alpha$ the assumptions of the theorem show

$$\varepsilon \geq 10\sqrt{\frac{\delta x}{n}} \geq \frac{60}{19}\sqrt{\frac{4\delta x}{n}} \geq \frac{60}{19}\sqrt{\frac{2(c\varepsilon^\alpha + \delta)x}{n}}.$$

Hence we have $\frac{19}{60}\varepsilon \geq \sqrt{\frac{2(c\varepsilon^\alpha+\delta)x}{n}}$ for all $\varepsilon$ satisfying the assumptions of the theorem. Now let $\tau := c\varepsilon^\alpha + \delta$ and $b := 2B$. By (35) and $\varepsilon \geq \frac{10Bx}{n}$ we then find

$$\frac{19\varepsilon}{20} > 3\mathbb{E}_{T'\sim P^n}\sup_{h\in\mathcal{H}}\mathbb{E}_{T'}h + \sqrt{\frac{2x\tau}{n}} + \frac{bx}{n}.$$

Applying Theorem 5.3 then yields

$$\mathrm{Pr}^*(T\in Z^n : \exists g\in\mathcal{G} \text{ with } \mathbb{E}_T g \leq \varepsilon/20 \text{ and } \mathbb{E}_P g \geq \varepsilon)$$

$$\leq P^n\bigg(T\in Z^n : \sup_{h\in\mathcal{H}}\mathbb{E}_T h \geq 19\varepsilon/20\bigg)$$

$$\leq P^n\bigg(T\in Z^n : \sup_{h\in\mathcal{H}}\mathbb{E}_T h > 3\mathbb{E}_{T'\sim P^n}\sup_{h\in\mathcal{H}}\mathbb{E}_{T'}h + \sqrt{\frac{2x\tau}{n}} + \frac{bx}{n}\bigg)$$

$$\leq e^{-x}. \qquad\qquad \square$$

With the help of the above lemma we can now prove the main result of this section, that is, Theorem 5.1.

PROOF OF THEOREM 5.1.   In order to apply Lemma 5.4 to the class $\mathcal{G}$ it obviously suffices to show the richness condition on $\mathcal{G}$ of Lemma 5.4. To this end let $f\in\mathcal{F}$ with

$$\mathbb{E}_T(L\circ f - L\circ f_{P,\mathcal{F}}) \leq \varepsilon/20 \quad\text{and}\quad \mathbb{E}_P(L\circ f - L\circ f_{P,\mathcal{F}}) \geq \varepsilon.$$

For $t\in[0,1]$ we define $f_t := tf + (1-t)f_{P,\mathcal{F}}$. Since $\mathcal{F}$ is convex we have $f_t\in\mathcal{F}$ for all $t\in[0,1]$. By the line-continuity of $L$ and Lebesgue's theorem we find that the map $h : t\mapsto \mathbb{E}_P(L\circ f_t - L\circ f_{P,\mathcal{F}})$ which maps from $[0,1]$ to $[0,B]$ is continuous. Since $h(0)=0$ and $h(1)\geq\varepsilon$ there is a $t\in(0,1]$ with

$$\mathbb{E}_P(L\circ f_t - L\circ f_{P,\mathcal{F}}) = h(t) = \varepsilon$$

by the intermediate value theorem. Moreover, for this $t$ we have

$$\mathbb{E}_T(L\circ f_t - L\circ f_{P,\mathcal{F}}) \leq \mathbb{E}_T(tL\circ f + (1-t)L\circ f_{P,\mathcal{F}} - L\circ f_{P,\mathcal{F}}) \leq \varepsilon/20.$$

Now, let $\varepsilon > 0$ with $\varepsilon \geq 10\max\{\omega_n(\mathcal{G}, c\varepsilon^\alpha + \delta), (\frac{\delta x}{n})^{1/2}, (\frac{4cx}{n})^{1/(2-\alpha)}, \frac{Bx}{n}\}$. Then by Lemma 5.4 we find that with probability at least $1 - e^{-x}$, every $f\in\mathcal{F}$ with $\mathbb{E}_T(L\circ f - L\circ f_{P,\mathcal{F}}) \leq \varepsilon/20$ satisfies $\mathbb{E}_P(L\circ f - L\circ f_{P,\mathcal{F}}) < \varepsilon$. Since we always have

$$\mathbb{E}_T(L\circ f_{T,\mathcal{F}} - L\circ f_{P,\mathcal{F}}) \leq 0 < \varepsilon/20,$$

we obtain the assertion.   $\square$

5.2. *Bounding the modulus of continuity.* The aim of this subsection is to bound the modulus of continuity of the class $\mathcal{G}$ in Theorem 5.1 with the help of covering numbers. We then present the resulting modification of Theorem 5.1.

Let us begin by recalling the definition of (local) Rademacher averages. To this end let $\mathcal{F}$ be a class of bounded measurable functions from $Z$ to $\mathbb{R}$ which is separable with respect to $\| \cdot \|_\infty$. Furthermore, let $P$ be a probability measure on $Z$ and $(\varepsilon_i)$ be a sequence of i.i.d. Rademacher variables (i.e., symmetric $\{-1, 1\}$-valued random variables) with respect to some probability measure $\mu$ on a set $\Omega$. Then the *Rademacher average* of $\mathcal{F}$ is

$$\mathrm{Rad}_P(\mathcal{F}, n) := \mathrm{Rad}(\mathcal{F}, n) := \mathbb{E}_{P^n} \mathbb{E}_\mu \sup_{f \in \mathcal{F}} \left| \frac{1}{n} \sum_{i=1}^n \varepsilon_i f(z_i) \right|,$$

and for $\varepsilon > 0$ the *local Rademacher average* of $\mathcal{F}$ is defined by

$$\mathrm{Rad}(\mathcal{F}, n, \varepsilon) := \mathrm{Rad}_P(\mathcal{F}, n, \varepsilon) := \mathbb{E}_{P^n} \mathbb{E}_\mu \sup_{\substack{f \in \mathcal{F}, \\ \mathbb{E}_P f^2 \le \varepsilon}} \left| \frac{1}{n} \sum_{i=1}^n \varepsilon_i f(z_i) \right|.$$

For a given $a > 0$ we immediately obtain $\mathrm{Rad}(a\mathcal{F}, n) = a \, \mathrm{Rad}(\mathcal{F}, n)$ and

$$(37) \qquad\qquad \mathrm{Rad}(a\mathcal{F}, n, \varepsilon) = a \, \mathrm{Rad}(\mathcal{F}, n, a^{-2}\varepsilon).$$

Moreover, by symmetrization the modulus of continuity can be estimated by the local Rademacher average. More precisely, we always have (see [36])

$$\omega_{P,n}(\mathcal{F}, \varepsilon) \le 2 \, \mathrm{Rad}_P(\mathcal{F}, n, \varepsilon), \qquad \varepsilon > 0.$$

Local Rademacher averages can be estimated by covering numbers. Without proof we state a slight modification of a corresponding result in [21]:

PROPOSITION 5.5. *Let $\mathcal{F}$ be a class of measurable functions from $Z$ to $[-1, 1]$ which is separable with respect to $\| \cdot \|_\infty$ and let $P$ be a probability measure on $Z$. Assume there are constants $a > 0$ and $0 < p < 2$ with*

$$\sup_{T \in Z^n} \log \mathcal{N}(\mathcal{F}, \varepsilon, L_2(T)) \le a\varepsilon^{-p}$$

*for all $\varepsilon > 0$. Then there exists a constant $c_p > 0$ depending only on $p$ such that for all $n \ge 1$ and all $\varepsilon > 0$ we have*

$$\mathrm{Rad}(\mathcal{F}, n, \varepsilon) \le c_p \max\left\{ \varepsilon^{1/2 - p/4} \left( \frac{a}{n} \right)^{1/2}, \left( \frac{a}{n} \right)^{2/(2+p)} \right\}.$$

Using this proposition we can replace the modulus of continuity in Theorem 5.1 by an assumption on the covering numbers of $\mathcal{G}$. Assuming that all resulting minimizers exist, the corresponding result then reads as follows:

THEOREM 5.6. *Let $\mathcal{F}$ be a convex set of bounded measurable functions from $Z$ to $\mathbb{R}$ and let $L : \mathcal{F} \times Z \to [0, \infty)$ be a convex and line-continuous loss function. For a probability measure $P$ on $Z$ we define*

$$\mathcal{G} := \{L \circ f - L \circ f_{P, \mathcal{F}} : f \in \mathcal{F}\}.$$

*Suppose that there are constants $c \geq 0$, $0 < \alpha \leq 1$, $\delta \geq 0$ and $B > 0$ with $\mathbb{E}_P g^2 \leq c(\mathbb{E}_P g)^\alpha + \delta$ and $\|g\|_\infty \leq B$ for all $g \in \mathcal{G}$. Furthermore, assume that $\mathcal{G}$ is separable with respect to $\|\cdot\|_\infty$ and that there are constants $a \geq 1$ and $0 < p < 2$ with*

$$\tag{38} \sup_{T \in Z^n} \log \mathcal{N}(B^{-1}\mathcal{G}, \varepsilon, L_2(T)) \leq a \varepsilon^{-p}$$

*for all $\varepsilon > 0$. Then there exists a constant $c_p > 0$ depending only on $p$ such that for all $n \geq 1$ and all $x \geq 1$ we have*

$$\mathrm{Pr}^*(T \in Z^n : \mathcal{R}_{L,P}(f_{T,\mathcal{F}}) > \mathcal{R}_{L,P}(f_{P,\mathcal{F}}) + c_p \varepsilon(n, a, B, c, \delta, x)) \leq e^{-x},$$

*where*

$$\varepsilon(n, a, B, c, \delta, x)$$
$$:= B^{2p/(4-2\alpha+\alpha p)} c^{(2-p)/(4-2\alpha+\alpha p)} \left(\frac{a}{n}\right)^{2/(4-2\alpha+\alpha p)} + B^{p/2} \delta^{(2-p)/4} \left(\frac{a}{n}\right)^{1/2}$$
$$+ B \left(\frac{a}{n}\right)^{2/(2+p)} + \sqrt{\frac{\delta x}{n}} + \left(\frac{cx}{n}\right)^{1/(2-\alpha)} + \frac{Bx}{n}.$$

PROOF. By (37) and Proposition 5.5 we find

$$\mathrm{Rad}(\mathcal{G}, n, \varepsilon) \leq c_p \max\left\{B^{p/2} \varepsilon^{1/2 - p/4} \left(\frac{a}{n}\right)^{1/2}, B\left(\frac{a}{n}\right)^{2/(2+p)}\right\}.$$

We assume without loss of generality that $c_p \geq 5$. Let $\varepsilon^* > 0$ be the largest real number that satisfies

$$\tag{39} \varepsilon^* = 2c_p B^{p/2} (c(\varepsilon^*)^\alpha + \delta)^{1/2 - p/4} \left(\frac{a}{n}\right)^{1/2}.$$

Furthermore, let $\varepsilon > 0$ be such that

$$\varepsilon = 2c_p \max\left\{B^{p/2}(c\varepsilon^\alpha + \delta)^{(2-p)/4} \left(\frac{a}{n}\right)^{1/2},\right.$$
$$\left. B\left(\frac{a}{n}\right)^{2/(2+p)}, \sqrt{\frac{\delta x}{n}}, \left(\frac{4cx}{n}\right)^{1/(2-\alpha)}, \frac{Bx}{n}\right\}.$$

It is easy to see that both $\varepsilon$ and $\varepsilon^*$ exist. Moreover, our above considerations show $\varepsilon \geq 10 \max\{\omega_n(\mathcal{G}, c\varepsilon^\alpha + \delta), (\frac{\delta x}{n})^{1/2}, (\frac{4cx}{n})^{1/(2-\alpha)}, \frac{Bx}{n}\}$, that is, $\varepsilon$ satisfies the assumptions of Theorem 5.1. In order to show the assertion it therefore suffices to bound $\varepsilon$ from above. To this end let us first assume that

$$B^{p/2}(c\varepsilon^\alpha + \delta)^{(2-p)/4}\left(\frac{a}{n}\right)^{1/2} \geq \max\left\{B\left(\frac{a}{n}\right)^{2/(2+p)}, \sqrt{\frac{\delta x}{n}}, \left(\frac{4cx}{n}\right)^{1/(2-\alpha)}, \frac{Bx}{n}\right\}.$$

Then we have $\varepsilon = 2c_pB^{p/2}(c\varepsilon^\alpha + \delta)^{(2-p)/4}(\frac{a}{n})^{1/2}$. Since $\varepsilon^*$ is the largest solution of this equation we hence find $\varepsilon \leq \varepsilon^*$. This shows that we always have

$$\varepsilon \leq \varepsilon^* + 2c_p\left(B\left(\frac{a}{n}\right)^{2/(2+p)} + \sqrt{\frac{\delta x}{n}} + \left(\frac{4cx}{n}\right)^{1/(2-\alpha)} + \frac{Bx}{n}\right).$$

Hence it suffices to bound $\varepsilon^*$ from above. To this end let us first assume $c(\varepsilon^*)^\alpha \geq \delta$. This implies $\varepsilon^* \leq 4c_pB^{p/2}(c \cdot (\varepsilon^*)^\alpha)^{1/2-p/4}(\frac{a}{n})^{1/2}$, and hence we find

$$\varepsilon^* \leq 16c_p^2B^{2p/(4-2\alpha+\alpha p)}c^{(2-p)/(4-2\alpha+\alpha p)}\left(\frac{a}{n}\right)^{2/(4-2\alpha+\alpha p)}.$$

Conversely, if $c(\varepsilon^*)^\alpha < \delta$ holds, then we immediately obtain

$$\varepsilon^* < 4c_pB^{p/2}\delta^{(2-p)/4}\left(\frac{a}{n}\right)^{1/2}. \qquad \square$$

**6. Variance bounds for SVMs.** In this section we prove some "variance bounds" in the sense of Theorem 5.6 for SVMs. Let us first ensure that these classifiers are ERM-type algorithms that fit into the framework of Theorem 5.6. To this end let $H$ be a RKHS of a continuous kernel over $X$, $\lambda > 0$, and $l: Y \times \mathbb{R} \to [0, \infty)$ be the hinge loss function. We define

$$(40) \qquad L(f, x, y) := \lambda\|f\|_H^2 + l(y, f(x))$$

and

$$(41) \qquad L(f, b, x, y) := \lambda\|f\|_H^2 + l(y, f(x) + b)$$

for all $f \in H$, $b \in \mathbb{R}$, $x \in X$ and $y \in Y$. Then $\mathcal{R}_{L,T}(\cdot)$ and $\mathcal{R}_{L,T}(\cdot, \cdot)$ obviously coincide with the objective functions of the SVM formulations and therefore SVMs are empirical $L$-risk minimizers. Furthermore note that all above minimizers exist (see [31]) and thus the SVM formulations in terms of $L$ actually fit into the framework of Theorem 5.6.

In the following, $f_{l,P}$ denotes a minimizer of $\mathcal{R}_{l,P}$ if no confusion can arise. For the shape of these minimizers which depend on $\eta := P(y = 1|\cdot)$ we refer to [39] and [30]. Now our first result is a variance bound which can be used when considering the empirical $l$-risk minimizer.

Lemma 6.1. *Let $P$ be a distribution on $X \times Y$ with Tsybakov noise exponent $0 \le q \le \infty$. Then there exists a minimizer $f_{l,P}$ mapping into $[-1, 1]$ such that for all bounded measurable functions $f : X \to \mathbb{R}$ we have*

$$\mathbb{E}_P(l \circ f - l \circ f_{l,P})^2$$

$$\le C_{\eta,q}(\|f\|_\infty + 1)^{(q+2)/(q+1)}(\mathbb{E}_P(l \circ f - l \circ f_{l,P}))^{q/(q+1)},$$

*where $C_{\eta,q} := \|(2\eta - 1)^{-1}\|_{q,\infty} + 2$ if $q > 0$ and $C_{\eta,q} = 1$ if $q = 0$.*

Proof. For $q = 0$ the assertion is trivial and hence we only consider the case $q > 0$. Given a fixed $x \in X$ we write $p := P(1|x)$ and $t := f(x)$. In addition, we introduce

$$v(p, t) := p(l(1, t) - l(1, f_{l,P}(x)))^2 + (1 - p)(l(-1, t) - l(-1, f_{l,P}(x)))^2,$$

$$m(p, t) := p(l(1, t) - l(1, f_{l,P}(x))) + (1 - p)(l(-1, t) - l(-1, f_{l,P}(x))).$$

Since Tsybakov's noise assumption implies $P_X(X_0) = 0$, we can restrict our consideration to $p \ne 1/2$. Now we will begin by showing

$$(42) \qquad v(p, t) \le \left(|t| + \frac{2}{|2p - 1|}\right) m(p, t).$$

Without loss of generality we may assume $p > 1/2$. Then we may set $f_{l,P}(x) := 1$ and thus we have $l(1, f_{l,P}(x)) = 0$ and $l(-1, f_{l,P}(x)) = 2$.

Let us first consider the case $t \in [-1, 1]$. Then we have $l(1, t) = 1 - t$ and $l(-1, t) = 1 + t$, and therefore (42) reduces to

$$(1 - t)^2 \le \left(|t| + \frac{2}{2p - 1}\right)(2p - 1)(1 - t).$$

Obviously, the latter inequality is equivalent to $1 - t \le (2p - 1)|t| + 2$, which is always satisfied for $t \in [-1, 1]$ and $p \ge 1/2$.

Now let us consider the case $t \le -1$. We then have $l(1, t) = 1 - t$ and $l(-1, t) = 0$, and after some elementary calculation we hence see that (42) is satisfied if and only if

$$p^2(6 - 2t) - p(5 - 3t) - 2t \ge 0.$$

The left-hand side is minimal if $p = (5 - 3t)/(12 - 4t)$, and thus we obtain

$$p^2(6 - 2t) - p(5 - 3t) - 2t \ge \frac{7t^2 - 18t - 25}{24 - 8t}.$$

Consequently, it suffices to show $7t^2 - 18t - 25 \ge 0$. However, the latter is true for all $t \le -1$ since $t \mapsto 7t^2 - 18t - 25$ is decreasing on $(-\infty, -1]$.

Now let us consider the third case, $t > 1$. Since we then have $l(1, t) = 0$ and $l(-1, t) = 1 + t$ it suffices to show

$$t - 1 \le t + \frac{2}{2p - 1}.$$

However, this is obviously true, and hence we have proved (42). Now, let us write

$$g(y, x) := l(y, f(x)) - l(y, f_{l,P}(x)),$$
$$h_1(x) := \eta(x)g(1, x) + (1 - \eta(x))g(-1, x),$$
$$h_2(x) := \eta(x)g^2(1, x) + (1 - \eta(x))g^2(-1, x).$$

Then (42) yields $h_2(x) \leq (\|f\|_\infty + \frac{2}{|2\eta(x)-1|})h_1(x)$ for all $x$ with $\eta(x) \neq 1/2$. For $t \geq 1$ we hence find

$$\mathbb{E}_P g^2 = \int_{|2\eta-1|^{-1} < t} h_2 \, dP_X + \int_{t \leq |2\eta-1|^{-1} < \infty} h_2 \, dP_X$$

$$\leq (\|f\|_\infty + 2t) \int_{|2\eta-1|^{-1} < t} h_1 \, dP_X + \int_{t \leq |2\eta-1|^{-1} < \infty} (\|f\|_\infty + 1)^2 \, dP_X$$

$$\leq 2(\|f\|_\infty + t)\mathbb{E}_P g + (\|f\|_\infty + 1)^2 P_X(|2\eta - 1|^{-1} \geq t)$$

$$\leq 2t(\|f\|_\infty + 1)\mathbb{E}_P g + (\|f\|_\infty + 1)^2 \|(2\eta - 1)^{-1}\|_{q,\infty} t^{-q}.$$

Let us define $t$ by $t^{q+1} := (\|f\|_\infty + 1)(\mathbb{E}_P g)^{-1}$. Since $\mathbb{E}_P g \leq \|f\|_\infty + 1$ we have $t \geq 1$ and hence the above estimate yields the assertion. $\square$

In the case of SVMs with offset we also need the following lemma which bounds the size of the offset $\tilde{b}_{P,\lambda}$. This lemma has been proved in [15] for empirical distributions. Although its generalization to general probability measures is straightforward we include the proof for completeness.

LEMMA 6.2. *Let $P$ be a distribution on $X \times Y$ and $\lambda > 0$. Then for all possible pairs $(\tilde{f}_{P,\lambda}, \tilde{b}_{P,\lambda}) \in H \times \mathbb{R}$ we have*

$$|\tilde{b}_{P,\lambda}| \leq \|\tilde{f}_{P,\lambda}\|_\infty + 1.$$

PROOF. If $P(y = y^*|x) = 1$ $P_X$-a.s. for some $y^* \in Y$, there is nothing to be proved since $\tilde{b}_{P,\lambda} = y^*$ by our assumption on SVMs mentioned in Section 2. Now let us assume that $\tilde{b}_{P,\lambda} > \|\tilde{f}_{P,\lambda}\|_\infty + 1$ and that $P$ is not degenerate in the above way. Then there exists a constant $\delta > 0$ such that $\tilde{b}_{P,\lambda} > \|\tilde{f}_{P,\lambda}\|_\infty + 1 + \delta$. This implies $\tilde{f}_{P,\lambda}(x) + \tilde{b}_{P,\lambda} > 1 + \delta$ for all $x \in X$. We define $b^*_{P,\lambda} := \tilde{b}_{P,\lambda} - \delta$. Obviously, we then find $l(-1, \tilde{f}_{P,\lambda}(x) + \tilde{b}_{P,\lambda}) = 0 = l(1, \tilde{f}_{P,\lambda}(x) + b^*_{P,\lambda})$ and

$$l(1, \tilde{f}_{P,\lambda}(x) + \tilde{b}_{P,\lambda}) = 1 + \tilde{f}_{P,\lambda}(x) + b^*_{P,\lambda} + \delta = l(-1, \tilde{f}_{P,\lambda}(x) + b^*_{P,\lambda}) + \delta$$

for all $x \in X$. Therefore we obtain $\mathcal{R}_{l,P}(\tilde{f}_{P,\lambda} + \tilde{b}_{P,\lambda}) > \mathcal{R}_{l,P}(\tilde{f}_{P,\lambda} + b^*_{P,\lambda})$ by using the assumption on $P$. $\square$

The proof of the above lemma can be easily generalized to a larger class of loss functions including, for example, the squared hinge loss. With the help of Lemma 6.1 we can now show a variance bound for SVMs. For brevity's sake we only state and prove the result for SVMs without offset. Therefore, the loss function $L$ is defined as in (40). Considering the proof, it is immediately clear that the variance bound also holds for the SVM with offset.

PROPOSITION 6.3. *Let $P$ be a distribution on $X \times Y$ with Tsybakov noise exponent $0 \le q \le \infty$. We define $C := 16 + 8\|(2\eta - 1)^{-1}\|_{q,\infty}$ if $q > 0$ and $C := 8$ otherwise. Furthermore, let $\lambda > 0$ and $0 < \gamma \le \lambda^{-1/2}$ such that $f_{P,\lambda} \in \gamma B_H$. Then for all $f \in \gamma B_H$ we have*

$$\mathbb{E}(L \circ f - L \circ f_{P,\lambda})^2 \le C(K\gamma + 1)^{(q+2)/(q+1)}(\mathbb{E}(L \circ f - L \circ f_{P,\lambda}))^{q/(q+1)}$$
$$+ 2C(K\gamma + 1)^{(q+2)/(q+1)}a^{q/(q+1)}(\lambda).$$

PROOF. We define $\hat{C} := (K\gamma + 1)^{(q+2)/(q+1)}$ and fix an $f \in \gamma B_H$. Furthermore, we choose a minimizer $f_{l,P}$ according to Lemma 6.1. Using $(a + b)^2 \le 2a^2 + 2b^2$ for all $a, b \in \mathbb{R}$ we first observe

$$\mathbb{E}(L \circ f - L \circ f_{P,\lambda})^2$$
$$\le 2\lambda^2\|f\|^4 + 2\lambda^2\|f_{P,\lambda}\|^4 + 2\mathbb{E}(l \circ f - l \circ f_{P,\lambda})^2$$
$$\le 4\mathbb{E}(l \circ f - l \circ f_{l,P})^2 + 4\mathbb{E}(l \circ f_{l,P} - l \circ f_{P,\lambda})^2 + 2\lambda^2\|f\|^4 + 2\lambda^2\|f_{P,\lambda}\|^4$$
$$\le 4C_{\eta,q}\hat{C}(\mathbb{E}(l \circ f - l \circ f_{l,P}) + \mathbb{E}(l \circ f_{P,\lambda} - l \circ f_{l,P}))^{q/(q+1)}$$
$$+ 2\lambda^2\|f\|^4 + 2\lambda^2\|f_{P,\lambda}\|^4,$$

where in the last step we have used Lemma 6.1 and $a^p + b^p \le 2(a + b)^p$ for all $a, b \ge 0$, $0 < p \le 1$. Since $\lambda\|f\|^2 \le 1$ and $\lambda\|f_{P,\lambda}\|^2 \le 1$, we can continue,

$$\mathbb{E}(L \circ f - L \circ f_{P,\lambda})^2$$
$$\le C\hat{C}\Big(\mathbb{E}(l \circ f - l \circ f_{l,P})$$
$$+ \mathbb{E}(l \circ f_{P,\lambda} - l \circ f_{l,P}) + \lambda^2\|f\|^4 + \lambda^2\|f_{P,\lambda}\|^4\Big)^{q/(q+1)}$$
$$\le C\hat{C}(\mathbb{E}(L \circ f - L \circ f_{P,\lambda}) + 2\mathbb{E}(l \circ f_{P,\lambda} - l \circ f_{l,P}) + 2\lambda\|f_{P,\lambda}\|^2)^{q/(q+1)}$$
$$\le C\hat{C}(\mathbb{E}(L \circ f - L \circ f_{P,\lambda}))^{q/(q+1)} + 2C\hat{C}a^{q/(q+1)}(\lambda). \qquad \square$$

**7. Proof of Theorem 2.8.** In this last section we prove our main result, Theorem 2.8. Since the proof is rather complex we split it into three parts. In Section 7.1 we estimate some covering numbers related to SVMs and Theorem 5.6. In Section 7.2 we then show that the trivial bound $\|f_{T,\lambda}\| \le$ $\lambda^{-1/2}$ can be significantly improved under the assumptions of Theorem 2.8. Finally, in Section 7.3 we prove Theorem 2.8.

7.1. *Covering numbers related to SVMs.* In this subsection we establish a simple lemma that estimates the covering numbers of the class $\mathcal{G}$ in Theorem 5.6 in terms of the covering numbers of $B_H$. For brevity's sake it only treats the case of SVMs with offset. The other case can be shown completely analogously.

LEMMA 7.1. *Let $H$ be a RKHS over $X$ such that $K$ defined by (12) satisfies $K \geq 1/2$, $P$ be a probability measure on $X \times Y$, $\lambda > 0$, and $L$ be defined by (41). Furthermore, let $1 \leq \gamma \leq \lambda^{-1/2}$ and*

$$\mathcal{F} := \{(f,b) \in H \times \mathbb{R} : \|f\|_H \leq \gamma \text{ and } |b| \leq \gamma K + 1\}.$$

*Defining $B := 2\gamma K + 3$ and $\mathcal{G} := \{L \circ (f,b) - L \circ (f_{P,\mathcal{F}}, b_{P,\mathcal{F}}) : (f,b) \in \mathcal{F}\}$ then gives $\|g\|_\infty \leq B$ for all $g \in \mathcal{G}$, where $(f_{P,\mathcal{F}}, b_{P,\mathcal{F}})$ denotes a $L$-risk minimizer in $\mathcal{F}$. Assume that there are constants $a \geq 1$ and $0 < p < 2$ such that for all $\varepsilon > 0$ we have*

$$\sup_{T \in Z^n} \log \mathcal{N}(B_H, \varepsilon, L_2(T_X)) \leq a\varepsilon^{-p}.$$

*Then there exists a constant $c_p > 0$ depending only on $p$ such that for all $\varepsilon > 0$ we have*

$$\sup_{T \in Z^n} \log \mathcal{N}(B^{-1}\mathcal{G}, \varepsilon, L_2(T)) \leq c_p a\varepsilon^{-p}.$$

PROOF. Let us write $\hat{\mathcal{G}} := \{L \circ (f,b) : (f,b) \in \mathcal{F}\}$ and $\mathcal{H} := \{l \circ (f + b) : (f,b) \in \mathcal{F}\}$. Furthermore, for brevity's sake we denote the set of all constant functions from $X$ to $[a,b]$ by $[a,b]$. We then have

$$\mathcal{N}(B^{-1}\mathcal{G}, \varepsilon, L_2(T)) = \mathcal{N}(B^{-1}\hat{\mathcal{G}}, \varepsilon, L_2(T)) \leq \mathcal{N}([0, \lambda\gamma^2] + B^{-1}\mathcal{H}, \varepsilon, L_2(T)).$$

Using the Lipschitz-continuity of the hinge loss and the subadditivity of the log-covering numbers we hence find

$$\log \mathcal{N}(B^{-1}\mathcal{G}, 3\varepsilon, L_2(T))$$
$$\leq \log \mathcal{N}([0, \lambda\gamma^2], \varepsilon, \mathbb{R}) + \log \mathcal{N}(B^{-1}\mathcal{H}, 2\varepsilon, L_2(T))$$
$$\leq \log\left(\frac{1}{\varepsilon} + 1\right) + \log \mathcal{N}(B^{-1}(\gamma \cdot B_H + [-B, B]), 2\varepsilon, L_2(T_X))$$
$$\leq 2\log\left(\frac{2}{\varepsilon} + 1\right) + \log \mathcal{N}(B_H, \varepsilon, L_2(T_X)).$$

From this we easily deduce the assertion. $\square$

7.2. *Shrinking the size of the SVM minimizers.* In this subsection we show that the trivial bound $\|f_{T,\lambda}\| \leq \lambda^{-1/2}$ can be significantly improved under the assumptions of Theorem 2.8. In view of Theorem 5.6 this improvement will have a substantial impact on the rates of Theorem 2.8. In order to obtain a rather flexible result let us suppose that for all $0 < p < 2$ we can determine constants $c, \gamma > 0$ such that

$$(43) \qquad \sup_{T \in Z^n} \log \mathcal{N}(B_{H_\sigma}, \varepsilon, L_2(T_X)) \leq c\sigma^{\gamma d}\varepsilon^{-p}$$

holds for all $\varepsilon > 0$, $\sigma \geq 1$. Recall that by Theorem 2.1 we can actually choose $\gamma := (1 - \frac{p}{2})(1 + \delta)$ for all $\delta > 0$.

LEMMA 7.2. *Let $X$ be the closed unit ball of the Euclidean space $\mathbb{R}^d$, and $P$ be a distribution on $X \times Y$ with Tsybakov noise exponent $0 \leq q \leq \infty$ and geometric noise exponent $0 < \alpha < \infty$. Furthermore, let us assume that (43) is satisfied for some $0 < \gamma \leq 2$ and $0 < p < 2$. Given an $0 \leq \varsigma < \frac{1}{5}$ we define*

$$\lambda_n := n^{-(4(\alpha+1)(q+1))/((2\alpha+1)(2q+pq+4)+4\gamma(q+1)) \cdot 1/(1-\varsigma)}$$

*and $\sigma_n := \lambda_n^{-1/((\alpha+1)d)}$. Assume that for the SVM without offset using the Gaussian RBF kernel with width $\sigma_n$ there are constants $\frac{1}{2(\alpha+1)} + 4\varsigma < \rho \leq \frac{1}{2}$ and $C \geq 1$ such that*

$$\mathrm{Pr}^*(T \in (X \times Y)^n : \|f_{T,\lambda_n}\| \leq Cx\lambda_n^{-\rho}) \geq 1 - e^{-x}$$

*for all $n \geq 1$ and all $x \geq 1$. Then there is another constant $\hat{C} \geq 1$ such that for $\hat{\rho} := \frac{1}{2}(\frac{1}{2(\alpha+1)} + 4\varsigma + \rho)$ and for all $n \geq 1$, $x \geq 1$ we have*

$$\mathrm{Pr}^*\left(T \in (X \times Y)^n : \|f_{T,\lambda_n}\| \leq \hat{C}x\lambda_n^{-\hat{\rho}}\right) \geq 1 - e^{-x}.$$

*Moreover, the same result is true for SVMs with offset.*

PROOF. We only prove the lemma for SVMs without offset since the proof for SVMs with offset is analogous. Now let $\hat{f}_{T,\lambda_n}$ be a minimizer of $\mathcal{R}_{L,T}$ on $Cx\lambda_n^{(\rho-1)/2}B_{H_{\sigma_n}}$, where $L$ is defined by (40). By our assumption we have $\hat{f}_{T,\lambda_n} = f_{T,\lambda_n}$ with probability not less than $1 - e^{-x}$ since $f_{T,\lambda_n}$ is unique for every training set $T$ by the strict convexity of $L$. We show that for some constant $\tilde{C} > 0$ and all $n \geq 1$, $x \geq 1$ the improved bound

$$(44) \qquad \|\hat{f}_{T,\lambda_n}\| \leq \tilde{C}x\lambda_n^{(\hat{\rho}-1)/2}$$

holds with probability not less than $1 - e^{-x}$. This then yields $\|f_{T,\lambda_n}\| \leq \tilde{C}x\lambda_n^{(\hat{\rho}-1)/2}$ with probability not less than $1 - 2e^{-x}$, and from the latter we

easily obtain the assertion. In order to establish (44) we will apply Theorem 5.6 to the modified SVM classifier which produces $\hat{f}_{T,\lambda_n}$. To this end we first remark that the infinite sample version $\hat{f}_{P,\lambda_n}$ which minimizes $\mathcal{R}_{L,P}$ on $Cx\lambda_n^{(\rho-1)/2}B_{H_{\sigma_n}}$ exists by a small modification of [31], Lemma 3.1. Furthermore, by Proposition 6.3 and assumption (43) we observe that we may choose $B$, $a$ and $c$ such that

$$B \sim x\lambda_n^{-\rho}, \qquad a \sim \lambda_n^{-\gamma/(\alpha+1)}, \qquad c \sim x^{(q+2)/(q+1)}\lambda_n^{-\rho\cdot(q+2)/(q+1)}.$$

In addition, Theorem 2.7 shows $a_{\sigma_n}(\lambda_n) \preceq \lambda_n^{\alpha/(\alpha+1)}$ and thus by Proposition 6.3 we may choose

$$\delta \sim x^{(q+2)/(q+1)}\lambda_n^{(\alpha q - \rho(q+2)(\alpha+1))/((\alpha+1)(q+1))}.$$

A rather time consuming but simple calculation then shows that the $\varepsilon$-term in Theorem 5.6 satisfies

$$\varepsilon(n,a,B,c,\delta,x) \preceq x^2 \lambda_n^{\frac{\alpha}{\alpha+1} - \frac{2\rho(\alpha+1)-1}{2(\alpha+1)} - \varsigma\frac{(2\alpha+1)(2q+pq+4)+4\gamma(q+1)}{2(\alpha+1)(2q+pq+4)}}.$$

Moreover, by Theorem 5.6 there is a constant $\tilde{C}_1 > 0$ independent of $n$ and $x$ such that for all $n \geq 1$ and all $x \geq 1$ the estimate

$$\lambda_n\|\hat{f}_{T,\lambda_n}\|^2 \leq \lambda_n\|\hat{f}_{T,\lambda_n}\|^2 + \mathcal{R}_{l,P}(\hat{f}_{T,\lambda_n}) - \mathcal{R}_{l,P}$$
$$\leq \lambda_n\|\hat{f}_{P,\lambda_n}\|^2 + \mathcal{R}_{l,P}(\hat{f}_{P,\lambda_n}) - \mathcal{R}_{l,P} + \tilde{C}_1 x^2\varepsilon(n,a,B,c,\delta,x)$$

holds with probability not less than $1 - e^{-x}$. Now $\lambda\|f_{P,\lambda}\|^2 \leq a_{\sigma_n}(\lambda_n) \preceq \lambda_n^{\alpha/(\alpha+1)}$ yields $\|f_{P,\lambda_n}\| \preceq \lambda_n^{-1/(2(\alpha+1))}$ and hence $\rho > \frac{1}{2(\alpha+1)}$ implies $\|f_{P,\lambda_n}\| \leq \lambda_n^{-\rho} \leq Cx\lambda_n^{-\rho}$ for large $n$. In other words, for large $n$ we have $f_{P,\lambda_n} = \hat{f}_{P,\lambda_n}$. Consequently, with probability not less than $1 - e^{-x}$ we have

$$\lambda_n\|\hat{f}_{T,\lambda_n}\|^2 \leq \lambda_n\|f_{P,\lambda_n}\|^2 + \mathcal{R}_{l,P}(f_{P,\lambda_n}) - \mathcal{R}_{l,P} + \tilde{C}_1 x^2\varepsilon(n,a,B,c,\delta,x)$$
$$\leq \tilde{C}_2 \lambda_n^{\alpha/(\alpha+1)} + \tilde{C}_1 x^2 \lambda_n^{\alpha/(\alpha+1)-(2\rho(\alpha+1)-1)/(2(\alpha+1))-4\varsigma},$$

which shows the assertion. $\quad\square$

7.3. *Proof of Theorem* 2.8. The next theorem almost establishes the result of Theorem 2.8. We present this intermediate result because it clarifies the impact of covering number bounds of the form (43) on our rates.

THEOREM 7.3. *Let $X$ be the closed unit ball of the Euclidean space $\mathbb{R}^d$, and $P$ be a distribution on $X \times Y$ with Tsybakov noise exponent $0 \leq q \leq \infty$ and geometric noise exponent $0 < \alpha < \infty$. Finally, let us assume that we can bound the covering numbers by (43) for some $0 < \gamma \leq 2$ and $0 < p < 2$. Given an $0 \leq \varsigma < \frac{1}{5}$ we define $\lambda_n$ and $\sigma_n$ as in Lemma 7.2. Then for all $\varepsilon > 0$ there* *is a constant $C > 0$ such that for all $x \geq 1$ and all $n \geq 1$ the SVM without offset and with regularization parameter $\lambda_n$ and Gaussian RBF kernel with width $\sigma_n$ satisfies*

$$\mathrm{Pr}^*(T : \mathcal{R}_P(f_{T,\lambda_n})$$
$$\leq \mathcal{R}_P + Cx^2 n^{-(4\alpha(q+1))/((2\alpha+1)(2q+pq+4)+4\gamma(q+1))\cdot 1/(1-\varsigma)+20\varsigma+\varepsilon)}$$
$$\geq 1 - e^{-x}.$$

*Moreover, the same result is true for SVMs with offset.*

PROOF.    Iteratively using Lemma 7.2 we find a constant $C \geq 1$ such that for $\rho := \frac{1}{2(\alpha+1)} + 4\varsigma + \varepsilon$ and all $n \geq 1$, $x \geq 1$ we have

$$\mathrm{Pr}^*(T \in (X \times Y)^n : \|f_{T,\lambda_n}\| \leq Cx\lambda_n^{-\rho}) \geq 1 - e^{-x}.$$

Repeating the calculations of Lemma 7.2 we hence find a constant $\tilde{C} > 0$ such that for all $n \geq 1$ and all $x \geq 1$ we have

$$\lambda_n \|f_{T,\lambda_n}\|^2 + \mathcal{R}_{l,P}(f_{T,\lambda_n}) - \mathcal{R}_{l,P}$$
$$\leq \lambda_n \|f_{P,\lambda_n}\|^2 + \mathcal{R}_{l,P}(f_{P,\lambda_n}) - \mathcal{R}_{l,P}$$
$$+ \tilde{C}_1 x^2 \lambda_n^{\alpha/(\alpha+1)-(2\rho(\alpha+1)-1)/(2(\alpha+1))-4\varsigma}$$

with probability not less than $1 - e^{-x}$. By the definition of $\rho$ we obtain

$$\lambda_n^{\alpha/(\alpha+1)-(2\rho(\alpha+1)-1)/(2(\alpha+1))-4\varsigma}$$
$$\leq \lambda_n^{\alpha/(\alpha+1)-4\varsigma-\varepsilon-4\varsigma}$$
$$\leq n^{-(4\alpha(q+1))/((2\alpha+1)(2q+pq+4)+4\gamma(q+1))\cdot 1/(1-\varsigma)+20\varsigma+3\varepsilon}.$$

From this we easily deduce the assertion.    $\square$

In order to prove Theorem 2.8 recall that by Theorem 2.1 we can choose $\gamma := (1 - \frac{p}{2})(1 + \delta)$ for all $\delta > 0$. The idea of the proof of Theorem 2.8 is to let $\delta \to 0$ while simultaneously adjusting $\varsigma$. The resulting rate is then optimized with respect to $p$. Unfortunately, a rigorous proof requires $p$ to be chosen a priori. Therefore, the optimization step is somewhat hidden in the following proof.

PROOF OF THEOREM 2.8.    Let us first consider the case $\alpha \leq \frac{q+2}{2q}$. Our aim is to apply Theorem 7.3. To this end we write $p_\delta := 2 - \delta$ and $\gamma_\delta := (1 - \frac{p_\delta}{2})(1 + \delta) = \frac{\delta}{2}(1 + \delta)$ for $\delta > 0$. Furthermore, we define $\varsigma_\delta$ by

$$\frac{4(\alpha+1)(q+1)}{(2\alpha+1)(4q-\delta q+4)+4\gamma_\delta(q+1)} \cdot \frac{1}{1-\varsigma_\delta} = \frac{\alpha+1}{2\alpha+1}.$$

Since $2\alpha q - q - 2 \leq 0 < 2\delta(q+1)$ we have $q(2\alpha+1) < 2(1+\delta)(q+1)$ and hence

$$4(2\alpha+1)(q+1) < 4(2\alpha+1)(q+1) - \delta q(2\alpha+1) + 2\delta(1+\delta)(q+1).$$

This shows $\varsigma_\delta > 0$ for all $\delta > 0$. Furthermore, these definitions also imply $\varsigma_\delta \to 0$ and $\gamma_\delta \to 0$ whenever $\delta \to 0$. Now Theorem 7.3 tells us that for all $\varepsilon > 0$ and all small enough $\delta > 0$ there exists a constant $C_{\delta,\varepsilon} \geq 1$ such that for all $n \geq 1$, $x \geq 1$ we have

$$\Pr{}^*(\mathcal{R}_P(f_{T,\lambda_n})$$
$$\leq \mathcal{R}_P + C_{\delta,\varepsilon}x^2 n^{-(4\alpha(q+1))/((2\alpha+1)(4q-\delta q+4)+4\gamma_\delta(q+1))\cdot 1/(1-\varsigma_\delta)+20\varsigma_\delta+\varepsilon)}$$
$$\geq 1 - e^{-x}.$$

In particular, if we choose $\delta$ sufficiently small we obtain the assertion.

Let us now consider the case $\frac{q+2}{2q} < \alpha < \infty$. In this case we write $p_\delta := \delta$ and $\gamma_\delta := (1 - \frac{p_\delta}{2})(1+\delta) = 1 + \frac{\delta}{2} - \frac{\delta^2}{2}$ for $\delta > 0$. Furthermore, we define $\varsigma_\delta$ by

$$\frac{4(\alpha+1)(q+1)}{(2\alpha+1)(2q+\delta q+4)+4\gamma_\delta(q+1)} \cdot \frac{1}{1-\varsigma_\delta} = \frac{2(\alpha+1)(q+1)}{2\alpha(q+2)+3q+4}.$$

Since for $0 < \delta \leq 1$ we have $0 < \delta q(2\alpha+1) + 2\delta(q+1) - 2\delta^2(q+1)$ we easily check that $\varsigma_\delta > 0$. Furthermore, the definitions ensure $\varsigma_\delta \to 0$ and $\gamma_\delta \to 1$ whenever $\delta \to 0$. The rest of the proof follows that of the first case. Finally, let us treat the case $\alpha = \infty$. We define $\alpha_\lambda$ by $\log \lambda = \alpha_\lambda d \log \frac{2\sqrt{d}}{\sigma}$. Since $\sigma > 2\sqrt{d}$ we have $\alpha_\lambda > 0$ for all $0 < \lambda < 1$. Furthermore, applying Theorem 2.7 for $\alpha_\lambda$ we find $a(\lambda) \leq 2C_d\lambda$ for all $0 < \lambda < 1$ and a constant $C_d > 0$ depending only on the dimension $d$. Adapted versions of Lemma 7.2 and Theorem 7.3 then yield the assertion. $\square$

**Acknowledgments.** We thank V. Koltchinskii and O. Bousquet for suggesting the local Rademacher averages as a way to obtain good performance bounds for SVMs and D. Hush for suggesting that "we are now in a position to obtain rates to Bayes." Finally we thank the anonymous reviewer for pointing out the example mentioned in Remark 2.12.

CCS-3
MS B256
LOS ALAMOS NATIONAL LABORATORY
LOS ALAMOS, NEW MEXICO 87545
USA
E-MAIL: ingo@lanl.gov

CCS-3
MS B265
LOS ALAMOS NATIONAL LABORATORY
LOS ALAMOS, NEW MEXICO 87545
USA
E-MAIL: jcs@lanl.gov